\documentclass[11pt]{amsart}

\usepackage{amssymb,amsmath,bm}
\usepackage{graphicx}
\usepackage{xcolor}
\usepackage{multirow}

\usepackage[top=26mm, bottom=26mm, left=32mm, right=32mm]{geometry}

      \usepackage{subfigure}

\newcommand{\vct}[1]{\bm{#1}}
\newcommand{\mtx}[1]{\mathsf{#1}}

\numberwithin{equation}{section}
\numberwithin{figure}{section}
\theoremstyle{definition}
\newtheorem{remark}{Remark}
\numberwithin{remark}{section}
\newtheorem{definition}{Definition}
\numberwithin{definition}{section}
\newcommand{\lsp}{\vspace{3mm}}

\newcommand{\vtwo}[2]{\left[\begin{array}{c} #1 \\ #2 \end{array}\right]}
\newcommand{\vthree}[3]{\left[\begin{array}{c} #1 \\ #2 \\ #3 \end{array}\right]}

\newcommand{\mtwo}[4]{\left[\begin{array}{cc}    #1 & #2 \\ #3 & #4  \end{array}\right]}
\newcommand{\mthree}[9]{\left[\begin{array}{ccc} #1 & #2 & #3 \\
                                                 #4 & #5 & #6 \\
                                                 #7 & #8 & #9 \end{array}\right]}

\newcommand{\yzcmt}[1]{\textcolor{black}{#1}}
\newcommand{\high}[1]{\textcolor{black}{#1}}

\begin{document}

\begin{center}

\textbf{A fast direct solver for two dimensional quasi-periodic multilayered  media scattering problems}

\lsp

\textit{\small Y. Zhang, and A. Gillman}

\lsp

\begin{minipage}{0.9\textwidth}\small
\noindent\textbf{Abstract:}
This manuscript presents a fast direct solution technique for solving two dimensional wave scattering problems 
from quasi-periodic multilayered structures. 
When the interface geometries are complex, 
the \yzcmt{dominant} term in \yzcmt{the} computational cost of creating the direct solver scales $O(NI)$ 
where $N$ is the number \yzcmt{of discretization} points on \high{each} interface and $I$ is the number of interfaces.  
The bulk of the  precomputation can be re-used for any choice of incident wave.  
As a result, the direct solver can solve over 200 scattering problems involving an eleven layer geometry 
with complex interfaces 100 times faster than building a new fast direct solver 
from scratch for each new set of boundary data. 
An added benefit of the presented solver is that building an \yzcmt{updated solver} for 
  a new geometry involving a replaced interface or \yzcmt{a} change in material property in \yzcmt{one} layer
  is inexpensive compared to building a new fast direct solver from scratch.
\end{minipage}
\end{center}
\section{Introduction}
\label{sec:intro}

\begin{figure}[htbp]
\centering
\begin{picture}(200,250)(100,0)
\put(0,0){\includegraphics[trim={0 2cm 0 2cm},clip,scale=0.7]{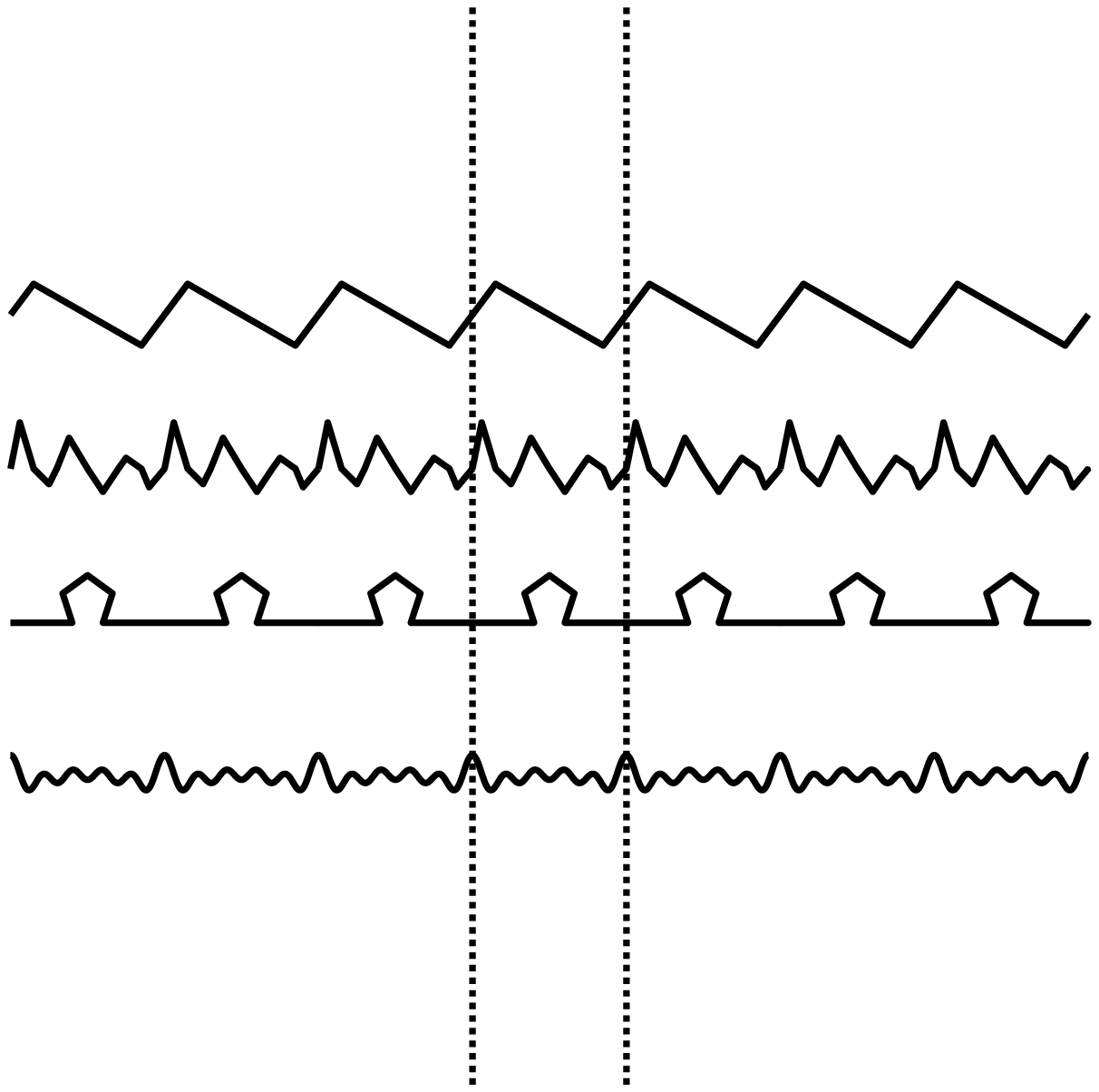}}
\put(100,220){\vector(1,-1){20}}
\put(100,210){\line(1,1){10}}
\put(105,205){\line(1,1){10}}
\put(110,230){$u^{inc}$}
\put(280,190){$\Omega_1$}
\put(280,145){$\Omega_2$}
\put(280,110){$\Omega_3$}
\put(280,75){$\Omega_4$}
\put(280,35){$\Omega_5$}
\put(330,160){$\Gamma_1$}
\put(330,130){$\Gamma_2$}
\put(330,90){$\Gamma_3$}
\put(330,60){$\Gamma_4$}
\put(196,30){\vector(-1,0){10}}
\put(210,30){\vector(1,0){10}}
\put(201,28){$d$}
\put(196,200){\vector(-1,0){10}}
\put(210,200){\vector(1,0){10}}
\put(201,198){$d$}
\end{picture}
\begin{center}
\caption{A five layered periodic geometry. 7 periods are shown.}
\label{fig:five_layer_geom}
\end{center}
\end{figure}

This manuscript considers the $I+1$ layered scattering problem defined by

\begin{equation}
 \label{eq:basic}
 \begin{split}
(\Delta + \omega_i^2) u_i(\vct{x}) &= 0  \  \qquad \vct{x}\in \Omega_i\\
 u_1-u_2 &= -u^{\rm inc}(\vct{x}) \  \qquad \vct{x}\in \Gamma_1\\
 \frac{\partial u_1}{\partial\nu} -\frac{\partial u_2}{\partial \nu} &= -\frac{\partial u^{\rm inc}}{\partial \nu} \  \qquad \vct{x}\in \Gamma_1\\
u_{i}-u_{i+1} &= 0 \  \qquad \vct{x}\in \Gamma_{i}, \ 1<i<I+1\\
 \frac{\partial u_i}{\partial\nu} -\frac{\partial u_{i+1}}{\partial \nu} &= 0 \ \qquad  \vct{x}\in \Gamma_{i},\ 1<i<I+1
 \end{split}
\end{equation}
\noindent
\yzcmt{where $u_i$ is the unknown solution in the region $\Omega_i\in\mathbb{R}^2$,
 the wave number in $\Omega_i$ is given by $\omega_i$ for $i = 1,\ldots,I+1$,
 and $\nu(\vct{x})$ is the normal vector at $\vct{x}$.}  
The interface $\Gamma_i$ for $i =1,\ldots,I$ between each layer is periodic
with period $d$.  
The boundary conditions enforce
continuity of the solution and its flux through the interfaces $\Gamma_i$.
The incident wave $u^{\rm inc}$ is defined by
$u^{\rm inc}(\vct{x}) = e^{{\rm i}\vct{k}\cdot\vct{x}}$
 where the incident vector is $\vct{k} = (\omega_1 \cos \theta^{\rm inc}, \omega_1 \sin \theta^{\rm inc})$
 and the incident angle is $-\pi<\theta^{\rm inc}<0$.   
Figure \ref{fig:five_layer_geom} illustrates a five layered periodic geometry.  The incident
wave $u^{\rm inc}$ is \emph{quasi-periodic} up to a phase, i.e. $u^{\rm inc}(x+d,y) = \alpha u^{\rm inc}(x,y)$
for $(x,y)\in \mathbb{R}^2$, where $\alpha$ is the Bloch phase defined by 
$$\alpha:= e^{i\omega_1 d \cos \theta^{\rm inc}}.$$
In the top and bottom layer, the solution \yzcmt{satisfies a radiation} 
 condition that is 
 characterized by the uniform convergence of the Rayleigh-Bloch expansions (see section
 \ref{sec:per} of \cite{BBStarling1994}.)

Multilayered periodic geometries are important in the design of 
optical and electromagnetic devices as well as select inverse scattering applications.
Some specific devices that involve scattering from \yzcmt{multilayered}
media are solar cells (thin-filmed photovoltaic cells \cite{Darbe2019SimulationAP,Miller,2010_solar1,2010_solar2} 
and solar 
thermal power \cite{Sergeant:10}), 
dielectric gratings for high-powered 
laser \cite{1995_diegrating,2004_laser,2019_dielectric} and wideband \cite{2010_chirped} 
applications.  
  Most of these applications 
require solving a scattering problem for a large number of incident 
angles $\theta^{\rm inc}$.  For example, in many engineering applications,
a Bragg diagram created from
the solution of 200 boundary value problems is desirable \cite{Barnett-Greengard}.  In optimal design
applications and inverse scattering, solving a scattering problem is nested inside of an optimization
loop.  At each step in the loop, a new scattering problem needs to be solved for many 
incident angles.  
When the geometry and material properties are close to the optimal choice or there are sufficient constraints on 
the material and/or interface properties the changes in the scattering problem are localized
to a few layer.

This paper presents a fast direct solver for the multilayered media
integral equation formulation presented in \cite{Cho:15}. This integral
formulation is robust even at \yzcmt{the} so-called Wood's anomalies.  The computational
cost of the proposed fast direct solver scales linearly with respect to the 
number of discretization points on the interfaces.  When a layer is changed
 with new material properties and/or a new interface \yzcmt{geometry}, the cost of updating the 
 direct solver scales linearly \yzcmt{with respect to} the number of discretization points
 affected by the modification.  For changing an interface, updating the 
 direct solver has a cost that scales linearly with respect to the number of 
 discretization points on the new interface.  For updating a \high{wave number} in $\Omega_i$,
 the cost scales linearly with respect to the number \yzcmt{of} discretization points on the 
 interfaces bounding $\Omega_i$.  This makes the solution technique a good option
 for \yzcmt{optimal design} and inverse scattering applications.


\subsection{Related work}

Direct discretization of (\ref{eq:basic}) is possible via finite difference or
finite element methods \cite{Hughes} but it faces two challenges\high{:} (i) meshing to the 
interfaces to maintain accuracy and (ii) enforcing the radiation condition.  Meshing \high{to the interfaces}  
can be effectively \yzcmt{handled} using mesh generation software such as GMSH \cite{gmesh}.  
Techniques such as perfectly matched layers \cite{2003_pml} can 
artificially enforce radiation conditions but can introduce artificial reflections and 
suffer from high condition numbers.  Another challenging aspect of using finite element methods
is that there is a loss in accuracy if the points per wavelength 
remains fixed (the so-called \emph{pollution} effect)
\cite{pollution}. 
 Another alternative method for directly discretizing (\ref{eq:basic})
is the rigorous-coupled wave analysis (RCWA) or Fourier Modal Method.  It is 
designed for multilayer gratings \cite{1981_RCWA} 
\yzcmt{and it depends on an iterative solve.} 
While a  
Fourier factorization method \cite{1996_RCWA_li1,1996_RCWA_li2} 
can be used to accelerate \yzcmt{the} convergence of an iterative solver, this solution 
approach is not ideal for problems with many right hand sides that arise in 
applications.  
\yzcmt{RCWA also has difficulty solving problems with interfaces that}
cannot be defined as a graph of a function
of the x-coordinate such as \yzcmt{the ``hedgehog'' interface in Figure \ref{fig:newinnterface}(b)}.  
Most often RCWA is used for geometries like cones and pillars. 
There is a concern that the method
is too simplified to capture complex structures \cite{RCWA1,RCWA2,RCWA3,RCWA4}.

When each layer is comprised of constant coefficient (i.e. not heterogeneous) medium,
it is possible to recast (\ref{eq:basic}) as a collection of boundary 
integral equations where the unknowns lie on the interfaces between
layers. 
There has been much work towards the use of boundary integral equations 
for quasi-periodic scattering problems including \cite{2019_Bruno,2002_Wilde,Arens_thesis,2006_Arens,2008_Nicholas}.
Reviews of boundary integral
equation techniques for scattering \yzcmt{off} a quasi-periodic array of obstacles are presented 
in \cite{Cho:15, Barnett-Greengard,2013_Gillman}.  
\high{The review in this paper focuses
on techniques for layered medium.} 
A boundary integral technique utilizing a fast direct solver for two layered media 
with periodic structures was presented in \cite{2014_periodic_HO}.  
 The integral formulation
utilized the quasi-periodic Green's function which is defined as an infinite
series.  For some choices of boundary data, this series does not converge
even though the problem is well posed.
  An incident angle $\theta^{\rm inc}$
that causes the quasi-periodic Green's function not to converge is called a \textit{Wood's anomaly}.  
There have been many techniques suggested to avoid these anomalies (such as \cite{Barnett-Greengard,2017_Bruno}).
{\color{black} A fast direct solver was constructed for quasi-periodic scattering off an 
infinite array of scatterers in \cite{2013_Gillman} for the robust integral formulation
presented in \cite{Barnett-Greengard}.  The work in this paper is an extension of that 
work to multilayered media problems.}  
This paper builds on the robust integral formulation in \cite{Cho:15} though it is
likely possible to build similar direct solvers for other formulations that are robust
at Wood's anomalies.  The integral formulation in \cite{Cho:15} makes
use of the free space Helmholtz Green's function, avoids the infinite sum and uses auxiliary unknowns to enforce
periodicity. {\color{black} The radiation condition is satisfied by enforcing continuity of the 
integral equation solution with the Rayleigh-Bloch expansions.}

Recently, \cite{Cho:18} replaced the boundary integral formulation in this 
approach to a \yzcmt{technique} based on the method of fundamental solutions.  This 
exchanges a second kind integral equation for a formulation 
that results in a system that is exponentially ill-conditioned.

\subsection{High level view of \yzcmt{the} solution technique}
Due to the problems associated with the quasi-periodic Green's function and {\color{black}a}
desire to exploit the constant coefficient medium, the 
fast direct solver is built for the robust boundary integral formulation 
proposed in \cite{Cho:15}.  Each interface has a boundary integral equation
that has ``structural'' similarities to a boundary integral equation for 
scattering off a {\color{black} single closed curve}.  The structural similarity is 
that a block matrix in the discretized boundary integral equation is amenable to fast direct inversion
techniques such as \textit{Hierarchically
Block Separable (HBS)} methods \cite{2012_martinsson_FDS_survey,2012_ho_greengard_fastdirect,2013_3DBIE}
which are closely related to the \textit{Hierarchically Semi-Separable (HSS)}
 \cite{2009_xia_superfast,2007_shiv_sheng,2010_xia}, the Hierarchical interpolative
 factorization (HIF) \cite{2014_HIF}, the $\mathcal{H}$ and $\mathcal{H}^{2}$-matrix methods
 \cite{2010_borm_book,2004_borm_hackbusch}.  
Roughly speaking these fast direct
 solvers utilize the fact that the off-diagonal blocks of the discretized integral 
 equation are low rank to create compressed representations of the matrix and 
 its inverse. 
 
 The linear system resulting from the discretization of the integral formulation
 in \cite{Cho:15} is rectangular where the principle sub-block is a block tridiagonal matrix.
 Each block in this tridiagonal matrix corresponds to a discretized boundary 
 integral operator {\color{black}that} (in the low frequency regime) is amenable to compression
 techniques such as those in fast direct solvers.  Utilizing this and separating 
 the matrices that depend on Bloch phase allows for the precomputation of the direct 
 solver to be utilized for all choices of incident angle.  The Bloch phase dependence 
 of many of the other block matrices in the rectangular
system can be separated out in a similar manner allowing them to be reused for 
multiple solves.  Further acceleration is gained by exploiting the block diagonal or
nearly block diagonal sparsity pattern of all the matrices. 
The combination of all these efforts dramatically reduces 
 the cost of processing the many solves needed in applications.

 {\color{black} 
The fast direct solver presented in this paper is ideally suited for applications that 
require many solves per geometry, involve solving problems where there are changes in 
a subset of the layers (material properties and/or interface geometries), or a combination
of many solves per geometry and changes in the geometry.  Applications where the solver
can be of benefit include optimal design of layered materials and inverse scattering problems
where the goal is \yzcmt{to} recover the thickness and/or the material properties of intermediate 
layers.  While the problems under consideration are acoustic scattering, the solution
technique can be extended to \textit{transverse electric} (TE) and \textit{transverse magnetic} 
(TM) wave problems.   

The direct solver presented in this paper is built for the  robust boundary integral formulation 
proposed in \cite{Cho:15} which enforces continuity of the solution and flux through 
interfaces.  The integral formulation can be extended to problems where there are jumps in
the solution and flux as long as these jumps are consistent
with the quasi-periodicity conditions.
 
 }

\subsection{Outline}
The paper begins by reviewing the integral formulation from \cite{Cho:15} in section 
\ref{sec:per}.  Next, section \ref{sec:FDS} presents the proposed fast direct solver.  
Numerical results in section \ref{sec:numerics} illustrate the performance of the direct 
solver.  Section \ref{sec:summary} summarizes manuscript and reviews the key features
of the presented work.

\section{Periodizing scheme}
\label{sec:per}

This section provides a review of the boundary integral formulation presented in 
\cite{Cho:15}.  The necessary integral operators are presented in \ref{sec:intop}.  Then
the full representation is presented in \ref{sec:intform}. Finally, the linear 
system resulting from enforcing continuity and quasi-periodicity of the solution 
is presented in section \ref{sec:sys}.

\begin{figure}[htbp]
\centering
\begin{picture}(200,250)(100,0)
\put(-50,0){\includegraphics[trim={0 0 0 0},clip,scale=0.5]{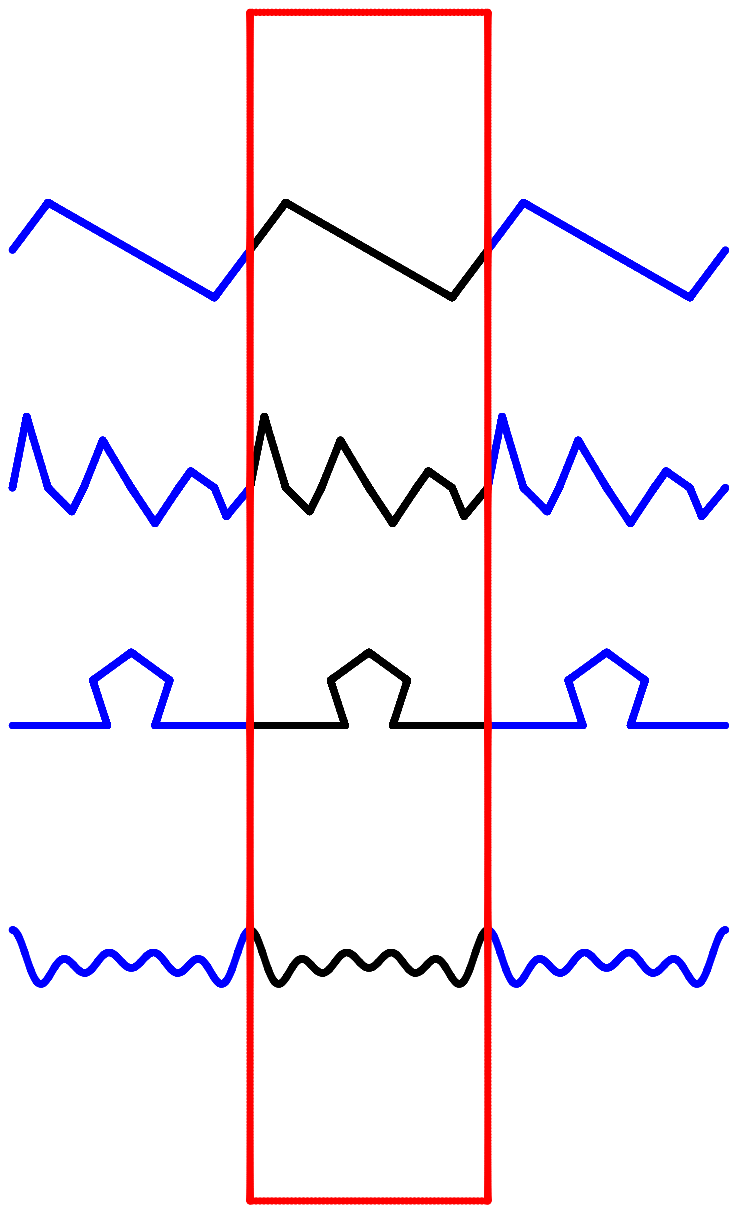}}
\put(90,-10){(a)}
\put(67,180){\color{red} $L$}
\put(118,180){\color{red} $R$}
\put(80,200){\color{red} $y=y_U$}
\put(80,10){\color{red} $y =y_D$}
\put(90,166){\footnotesize$\Gamma_1$}
\put(90,135){\footnotesize$\Gamma_2$}
\put(90,105){\footnotesize$\Gamma_3$}
\put(90,63){\footnotesize$\Gamma_4$}
\put(180,-30){\includegraphics[trim={5cm 0 5cm 0},clip,scale=0.65]{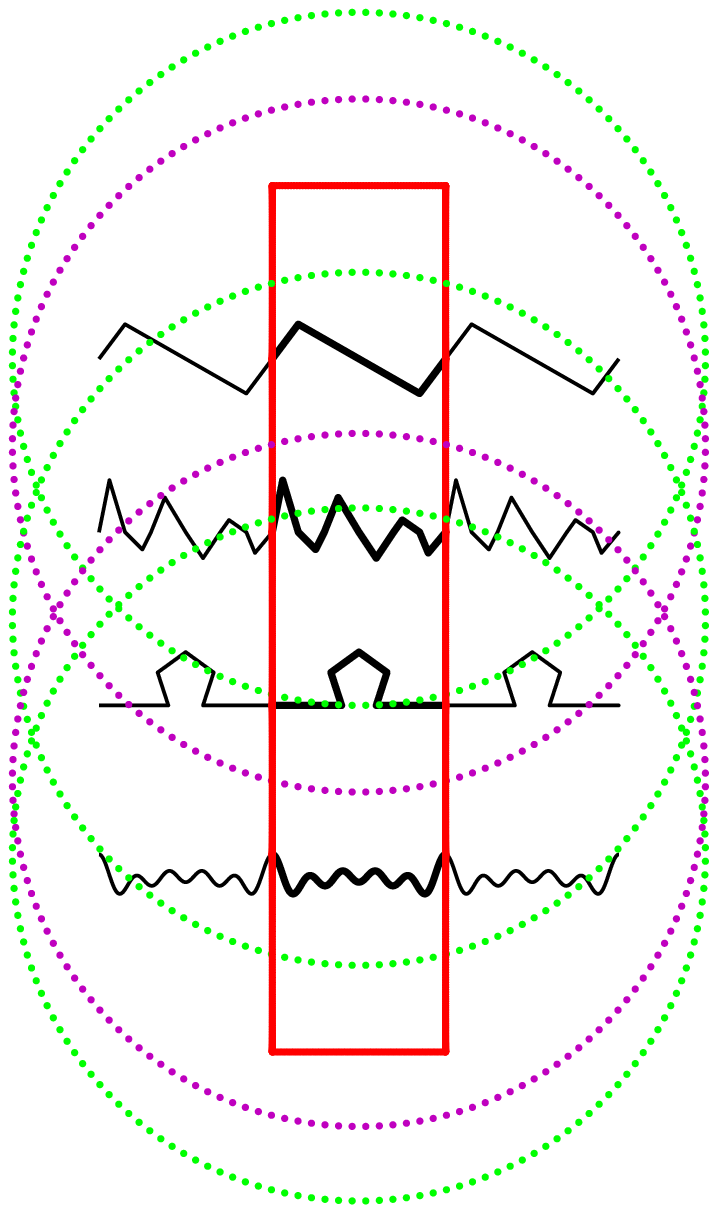}}
\put(280,-10){(b)}
\put(242,170){\footnotesize$\Omega_1$}
\put(242,140){\footnotesize$\Omega_2$}
\put(242,110){\footnotesize$\Omega_3$}
\put(242,75){\footnotesize$\Omega_4$}
\put(242,38){\footnotesize$\Omega_5$}
\put(300,222){\color{green}$P_1$}
\put(310,182){\color{magenta}$P_2$}
\put(345,112){\color{green}$P_3$}
\put(310,35){\color{magenta}$P_4$}
\put(303,-2){\color{green}$P_5$}

\end{picture}
\begin{center}
\caption{
This figure illustrates a five layered periodic geometry with artificial walls and proxy circles. 
Only three periods of the infinite periodic geometry are shown.  The period contained within the
unit cell is in black while the other two periods \yzcmt{are in blue}.
Figure (a) illustrates the notation for the unit cell with left, right, upper, and lower boundary $L$, $R$, $U$, and $D$ shown in red lines.
Figure (b) illustrates the proxy circles $P_i$ for each layer.  The color of the proxy circles alternates 
between green and magenta.
}
\label{fig:five_layer_geom_unitcell}
\end{center}
\end{figure}

The integral formulation proposed in \cite{Cho:15} solves (\ref{eq:basic}) in an infinite vertical unit strip
of width $d$.  Because the solution is known to be quasi-periodic, the solution outside of the 
unit strip can be found by scaling the solution by the appropriate Bloch phase factor. Let $x = L$ and $x = R$ 
denote the left and right bounds for the unit strip.  The solution technique further partitions space 
by introducing artificial top and bottom \yzcmt{walls} to the unit strip at $y = y_U$ and $y = y_D$ respectively.
Figure \ref{fig:five_layer_geom_unitcell}(a) illustrates this partitioning.  The box bounded by these artificial 
boundaries is called the \textit{unit cell}.
Inside the unit cell the solution is represented via an integral formulation.  Above and below the unit cell,
(i.e. for points in the unit strip where $y>y_U$ \yzcmt{or} $y<y_D$), the solution is given by \yzcmt{Rayleigh-Bloch} expansions.
Specifically, for $\vct{x} = (x,y)$ in the unit strip where $y>y_U$, the solution is given by 

\begin{equation}u(x,y) = \sum_{n\in\mathbb{Z}} a^U_n e^{{\rm i}\kappa_nx}e^{{\rm i}k^U_n(y-y_U)}
 \label{eq:ral1}
\end{equation}

and, for $\vct{x}=(x,y)$ in the unit strip where $y<y_D$, the solution is given by 
\begin{equation}u(x,y) = \sum_{n\in\mathbb{Z}} a_n^D e^{{\rm i}\kappa_nx}e^{{\rm i}k^D_n(-y+y_D)}
 \label{eq:ral2}
\end{equation}

where $\kappa_n:= \omega_1\cos\theta^{\rm inc} +\frac{2\pi n}{d}$, $k^U_n = \sqrt{\omega_1^2-\kappa_n^2}$,
$k^D_n = \sqrt{\omega_{I+1}^2-\kappa_n^2}$
and the sets $\{a^U_n\}$ and $\{a^D_n\}$ are coefficients to be determined.
\textcolor{black}{The square root can be either a positive real or positive imaginary number. }

 \subsection{Integral operators} 
 \label{sec:intop}
This section presents the integral operators needed to represent the solution inside
the unit cell.

Let $\Gamma_i$ for $i=1,\ldots, I$ denote the interfaces
inside the unit cell and $\Omega_i$ denote the regions in between
\yzcmt{for each layer} 
in the unit cell.  Both are numbered from the top down.  
Figure \ref{fig:five_layer_geom_unitcell}(a) illustrates the numbering 
of the five layered geometry within the unit cell.


Let $G_{\omega}(\vct{x},\vct{y}) = \frac{i}{4}H_0^{(1)}(\omega\|\vct{x}-\vct{y}\|)$ denote
the two dimensional free space Green's function for the Helmholtz equation with wave number 
$\omega$ where $H_0^1$ is the Hankel function of zeroth order \cite{Abramowitz-Stegun}.

The standard Helmholtz single and double layer integral operators defined on a curve $W$ \cite{coltonkress} are 
$$(\mathcal{S}_W^\omega\rho)(\vct{x}) = \int_W G_\omega(\vct{x},\vct{y})\rho(\vct{y})dl(\vct{y}) \quad 
\yzcmt{\mbox{and} \quad}
(\mathcal{D}_W^\omega\rho)(\vct{x}) = \int_W \partial_{\nu_{\vct{y}}} G_\omega(\vct{x},\vct{y})\rho(\vct{y})dl(\vct{y}),$$
respectively, where $\nu_{\vct{y}}$ denotes the normal vector at the point $\vct{y}\in W$.

For the periodizing scheme, integral operators involving the unit cell and its 
neighbors (left and right) are needed.  
 
These operators, denoted with tilde, are defined as follows
\begin{equation}
\begin{split}
 \label{eq:single} (\tilde{\mathcal{S}}_W^\omega\rho)(\vct{x}) &= \sum_{l=-1}^1 \alpha^{l}\int_W G_\omega(\vct{x},\vct{y}+l\vct{d})\rho(\vct{y})dl(\vct{y})\\
 &= (\mathcal{S}_W^\omega \rho)(\vct{x}) + {\color{black}(\mathcal{S}_{W}^{\omega,pm} \rho)}(\vct{x})
\end{split}
 \end{equation}
and 
\begin{equation}
 \label{eq:double}
 \begin{split}(\tilde{\mathcal{D}}_W^\omega\rho)(\vct{x}) &= \sum_{l=-1}^1 \alpha^{l}\int_W \partial_{\nu_{\vct{y}}} G_\omega(\vct{x},\vct{y}+l\vct{d})\rho(\vct{y})dl(\vct{y})\\
  &= (\mathcal{D}_W^\omega \rho)(\vct{x}) + ({\color{black}\mathcal{D}_{W}^{\omega,pm} \rho})(\vct{x})
 \end{split}
\end{equation}
where 
{\color{black}\begin{equation}
\label{eq:sing_pm} 
\displaystyle({\color{black}\mathcal{S}_{W}^{\omega,pm}\rho)}(\vct{x}) = \sum_{l=-1,l\neq 0}^1 \alpha^{l}\int_W G_\omega(\vct{x},\vct{y}+l\vct{d})\rho(\vct{y})dl(\vct{y})
\end{equation}}
and 
{\color{black}\begin{equation}
 \label{eq:doub_pm}
\displaystyle({\color{black}\mathcal{D}^{\omega,pm}_{W}}\rho)(\vct{x}) = \sum_{l=-1,l\neq 0}^1 \alpha^{l}\int_W \partial_{\nu_{\vct{y}}} G_\omega(\vct{x},\vct{y}+l\vct{d})\rho(\vct{y})dl(\vct{y}).
\end{equation}}
The notation ${pm}$ stands for plus-minus referring to the left and right neighboring copies of \yzcmt{$W$}.

These integral operators are not sufficient to enforce quasi-periodicity.  
They are missing information from the infinite copies that are ``far'' from the 
unit cell.  A proxy basis is used to capture the missing information.  For 
simplicity, consider a layer $\Omega_l$.  Let $\{\vct{y}_j\}_{j=1}^P$ denote
a collection of uniformly distributed points on a circle $P_l$ of radius $2d$ that is centered 
in $\Omega_l$.  {\color{black} The proxy circle needs to be large enough to shield the interface in
the unit cell \yzcmt{from} its far-field copies, 
which are more than $\frac{3d}{2}$ away from the center of $\Gamma_i$ in the horizontal direction. 
It is proved in \cite{timo} that larger proxy radius leads to higher order convergence rate with respect to the number of
basis functions $P$. However, the radius cannot be arbitrarily large, as the {\color{black} magnitude of the coefficients}
grows exponentially with respect to the ratio between the proxy radius and $\frac{3d}{2}$.
We set the radius of the proxy circle to be $R_{proxy}\in\left[ \frac{3d}{2}  , 2d\right]$ 
as in \cite{Cho:15}.}
The elements of the \textit{proxy basis} used to capture \yzcmt{the} far 
field information are defined by
\begin{equation}
 \phi^{\omega_l}_j = \frac{\partial G_{\omega_l}}{\partial \vct{n}_j} (\vct{x},\vct{y}_j) + i\omega_l G_{\omega_l}(\vct{x},\vct{y}_j)
 \label{eq:proxy}
\end{equation}
where $\vct{n}_j$ is the normal vector at $\vct{y}_j$ on $P_l$.  
This choice of basis results in smaller \yzcmt{coefficients} when compared to \yzcmt{using} just 
the single or double layer potential as a basis \cite{Cho:15}.  
If the layer has a high aspect ratio (i.e. taller than $d$), the proxy surface should be taken to be an 
ellipse; see page 8 of \cite{Cho:15}.
 Figure \ref{fig:five_layer_geom_unitcell}(b) illustrates the proxy circles for a five layered geometry.

The boundary integral equations involve additional integral operators which we define in 
this section for simplicity of presentation.  Specifically, an 
integral operator defined on an interface $W$ will need to be evaluated at $\vct{x}\in V$
where $V$ is an interface (the same or a vertical neighbor \yzcmt{of $W$}).  
For $\vct{x}\in V$ where $V$ is an interface, let $(\tilde{{S}}_{V,W}^\omega\rho)$ 
denote the evaluation of (\ref{eq:single}) at $\vct{x}$, i.e.
$$ (\tilde{{S}}_{V,W}^\omega\rho)(\vct{x}) = \sum_{l=-1}^1 \alpha^{l}\int_W G_\omega(\vct{x},\vct{y}+l\vct{d})\rho(\vct{y})dl(\vct{y}).$$
Likewise, let $(\tilde{{D}}_{V,W}^\omega\rho)$ 
denote the evaluation of \yzcmt{(\ref{eq:double})} at $\vct{x}\in V$, i.e.
$$ (\tilde{{D}}_{V,W}^\omega\rho)(\vct{x}) = \sum_{l=-1}^1 \alpha^{l}\int_W \partial_{\nu_{\vct{y}}} G_\omega(\vct{x},\vct{y}+l\vct{d})\rho(\vct{y})dl(\vct{y}).$$
The $pm$ notation for the neighbor interactions follows in a similar \yzcmt{fashion.} 
\yzcmt{For example, the operator} $(\tilde{{S}}_{V,W}^\omega\rho)(\vct{x}) $ can be written as the following sum
$$
(\tilde{{S}}_{V,W}^\omega\rho)(\vct{x})= ({S}_{V,W}^\omega\rho)(\vct{x}) + {\color{black}({S}_{V, W}^{\omega,pm}\rho)}(\vct{x}),
$$
where
$$
 ({\color{black}{S}_{V, W}^{\omega,pm}}\rho)(\vct{x})= \sum_{l=-1,\, l\neq 0}^1 \alpha^{l}\int_W  G_\omega(\vct{x},\vct{y}+l\vct{d})\rho(\vct{y})dl(\vct{y}).
$$

In order to enforce continuity of the fluxes, 
the normal derivatives of these integral operators are 
required.   For $\vct{x}\in V$ where $V$ is an interface, let $(\tilde{D}^{*,\omega}_{W,V}\rho)$
denote the evaluation of the normal derivative of the single layer operator (\ref{eq:single}) at $\vct{x}$, i.e.
$$(\tilde{{D}}_{W,V}^{*,\omega}\rho)(\vct{x}) = \sum_{l=-1}^1 \alpha^{l}\int_W \partial_{\nu_{\vct{x}}}G_\omega(\vct{x},\vct{y}+l\vct{d})\rho(\vct{y})dl(\vct{y})$$
where $\nu_{\vct{x}}$ is the normal vector at $\vct{x}\in V$.
Similarly, let $(\tilde{T}^\omega_{W,V}\rho)$ denote the evaluation of the normal derivative of 
the double layer operator (\ref{eq:double}) at $\vct{x}$, i.e.
$$(\tilde{{T}}_{W,V}^{\omega}\rho)(\vct{x}) = \sum_{l=-1}^1 \alpha^{l}\int_W \partial_{\nu_{\vct{x}}}\partial_{\nu_{\vct{y}}}G_\omega(\vct{x},\vct{y}+l\vct{d})\rho(\vct{y})dl(\vct{y}).$$ 

\subsection{Integral formulation}
\label{sec:intform}

The periodizing scheme within the unit cell is based on a modified version of the 
combined field boundary integral formulation \cite{1978_kress,rokh83}.  
Specifically, the solution in the unit cell is expressed as
\begin{align}
\label{eq:lev1}
   u_1(\vct{x}) &= (\tilde{\mathcal{S}}^{\omega_1}_{\Gamma_1}\sigma_1)(\vct{x}) + (\tilde{\mathcal{D}}^{\omega_1}_{\Gamma_1}\tau_1)(\vct{x})
 + \sum_{j=1}^Pc^1_j\phi^{\omega_1}_j(\vct{x}) \quad {\rm for} \ \vct{x}\in \Omega_1,\\
 \label{eq:levI} u_{I+1}(\vct{x}) &= (\tilde{\mathcal{S}}^{\omega_{I+1}}_{\Gamma_I}\sigma_I)(\vct{x}) + 
    (\tilde{\mathcal{D}}^{\omega_{I+1}}_{\Gamma_I}\tau_I)(\vct{x})
 + \sum_{j=1}^Pc^{I+1}_j\phi^{\omega_{I+1}}_j(\vct{x}) \ {\rm for} \ \vct{x}\in \Omega_{I+1}, \quad \high{\rm and}\\
 \label{eq:levmid}
u_i(\vct{x}) &= (\tilde{\mathcal{S}}^{\omega_i}_{\Gamma_{i-1}}\sigma_{i-1})(\vct{x}) + 
(\tilde{\mathcal{D}}^{\omega_i}_{\Gamma_{i-1}}\tau_{i-1})(\vct{x}) + 
(\tilde{\mathcal{S}}^{\omega_i}_{\Gamma_{i}}\sigma_{i})(\vct{x}) + 
(\tilde{\mathcal{D}}^{\omega_i}_{\Gamma_{i}}\tau_i)(\vct{x}) + 
\sum_{j=1}^Pc^i_j\phi^{\omega_i}_j(\vct{x})  
\end{align}
for $\vct{x}\in\Omega_i$, $2\leq i\leq I$ where $\sigma_i$ and $\tau_i$ are unknown boundary charge distributions
and $\{c^i_j\}_ {j=1}^P$ are unknown constants , for $i = 1,\dots,I$.

{\color{black}Enforcing the transmission condition in equation (\ref{eq:basic}) corresponding
to continuity of the solution through the interfaces results in the following \yzcmt{integral
equations:}}

\begin{equation}
\begin{split}
 \label{eq:DIEtop} 
-\tau_1 + (\tilde{D}^{\omega_1}_{\Gamma_1,\Gamma_1}&  - \tilde{D}^{\omega_2}_{\Gamma_1,\Gamma_1})\tau_1
 + (\tilde{S}^{\omega_1}_{\Gamma_1,\Gamma_1}  - \tilde{S}^{\omega_2}_{\Gamma_1,\Gamma_1})\sigma_1
 - \tilde{D}^{\omega_2}_{\Gamma_1,\Gamma_2} \tau_2 - \tilde{S}^{\omega_2}_{\Gamma_1,\Gamma_2}\sigma_2
  \\
& +\sum_{p = 1}^P(c^1_p \phi^{\omega_1}_p-c^2_p\phi^{\omega_2}_p)|_{\Gamma_1}= -u^{\rm inc} 
\ {\rm on} \ \Gamma_1,
\end{split}
\end{equation}

\begin{equation}
 \label{eq:DIEbot}
 \begin{split}-\tau_I + (\tilde{D}^{\omega_{I}}_{\Gamma_I,\Gamma_I}  &- \tilde{D}^{\omega_{I+1}}_{\Gamma_{I},\Gamma_I})\tau_I
 + (\tilde{S}^{\omega_I}_{\Gamma_I,\Gamma_I}  - \tilde{S}^{\omega_I}_{\Gamma_I,\Gamma_I})\sigma_I
  - \tilde{D}^{\omega_{I+1}}_{\Gamma_{I},\Gamma_{I-1}} \tau_{I-1} - \tilde{S}^{\omega_{I-1}}_{\Gamma_{I},\Gamma_{I-1}}\sigma_{I-1}\\
&+\sum_{p = 1}^P(c^I_p \phi^{\omega_I}_p -c^{I+1}_p\phi^{\omega_{I+1}}_p)|_{\Gamma_I} = 0
\ {\rm on} \ \Gamma_I,
\end{split}
\end{equation}

and 

\begin{equation}
 \label{eq:DIE}  
 \begin{split}- \tau_i +  (\tilde{D}^{\omega_{i}}_{\Gamma_i,\Gamma_i}  &- \tilde{D}^{\omega_{i+1}}_{\Gamma_{i},\Gamma_i})\tau_i
  + (\tilde{S}^{\omega_i}_{\Gamma_i,\Gamma_i}  - \tilde{S}^{\omega_i}_{\Gamma_i,\Gamma_i})\sigma_i 
  + \tilde{D}^{\omega_{i}}_{\Gamma_i,\Gamma_{i-1}}\tau_{i-1} + \tilde{D}^{\omega_{i+1}}_{\Gamma_i,\Gamma_{i+1}}\tau_{i+1}+
 \\
&  \tilde{S}^{\omega_{i}}_{\Gamma_i,\Gamma_{i-1}}\sigma_{i-1} + \tilde{S}^{\omega_{i+1}}_{\Gamma_i,\Gamma_{i+1}}\sigma_{i+1}+\sum_{p = 1}^P(c^i_p \phi^{\omega_i}_p -c^{i+1}_p\phi^{\omega_{i+1}}_p)|_{\Gamma_i} = 0
\ {\rm on} \ \Gamma_i \ { \rm for} \ 1< i < I
\end{split}
\end{equation}
where $\tilde{S}^{\omega_i}_{\Gamma_{i+1},\Gamma_{i}}$ denotes the 
periodized single layer integral operator (\ref{eq:single}) 
defined on $\Gamma_i$ evaluated on $\Gamma_{i+1}$, etc.

Likewise, enforcing the {\color{black} transmission condition in equation (\ref{eq:basic}) 
corresponding to continuity of the flux through}
 the interfaces results in the following boundary integral equations:

\begin{equation}
 \label{eq:NIEtop} \begin{split}-\sigma_1 + &(\tilde{T}^{\omega_1}_{\Gamma_1,\Gamma_1}  - \tilde{T}^{\omega_2}_{\Gamma_1,\Gamma_1})\tau_1
 + (\tilde{D}^{*,\omega_1}_{\Gamma_1,\Gamma_1}  - \tilde{D}^{*,\omega_2}_{\Gamma_1,\Gamma_1})\sigma_1
 - \tilde{T}^{\omega_2}_{\Gamma_1,\Gamma_2} \tau_2 - \tilde{D}^{*,\omega_2}_{\Gamma_1,\Gamma_2}\sigma_2\\
&+\sum_{p = 1}^P\left(c^1_p \frac{\partial\phi^{\omega_1}_p}{\partial \nu} -c^2_p\frac{\partial\phi^{\omega_2}_p}{\partial\nu}\right)|_{\Gamma_1}
= -u^{\rm inc} 
\ {\rm on} \ \Gamma_1,
\end{split}
\end{equation}

\begin{equation}
 \label{eq:NIEbot} \begin{split}-\sigma_I &+ (\tilde{T}^{\omega_{I}}_{\Gamma_I,\Gamma_I}  - \tilde{T}^{\omega_{I+1}}_{\Gamma_{I},\Gamma_I})\tau_I
 + (\tilde{D}^{*,\omega_I}_{\Gamma_I,\Gamma_I}  - \tilde{D}^{*,\omega_I}_{\Gamma_I,\Gamma_I})\sigma_I
  - \tilde{T}^{\omega_{I+1}}_{\Gamma_{I},\Gamma_{I-1}} \tau_{I-1} - \tilde{D}^{*\omega_{I-1}}_{\Gamma_{I},\Gamma_{I-1}}\sigma_{I-1}\\
&+\sum_{p = 1}^P\left(c^I_p \frac{\partial\phi^{\omega_I}_p}{\partial \nu} -c^{I+1}_p\frac{\partial\phi^{\omega_{I+1}}_p}{\partial \nu}\right)|_{\Gamma_I} = 0
\ {\rm on} \ \Gamma_I,
\end{split}
\end{equation}

and 

\begin{equation}
 \label{eq:NIE}  \begin{split}- \sigma_i &+  (\tilde{T}^{\omega_{i}}_{\Gamma_i,\Gamma_i}  - \tilde{T}^{\omega_{i+1}}_{\Gamma_{i},\Gamma_i})\tau_i
  + (\tilde{D}^{*,\omega_i}_{\Gamma_i,\Gamma_i}  - \tilde{D}^{*,\omega_i}_{\Gamma_i,\Gamma_i})\sigma_i + 
  + \tilde{T}^{\omega_{i}}_{\Gamma_i,\Gamma_{i-1}}\tau_{i-1} + \tilde{T}^{\omega_{i+1}}_{\Gamma_i,\Gamma_{i+1}}\tau_{i+1} \\
  &+ \tilde{D}^{*,\omega_{i}}_{\Gamma_i,\Gamma_{i-1}}\sigma_{i-1} + \tilde{D}^{*,\omega_{i+1}}_{\Gamma_i,\Gamma_{i+1}}\sigma_{i+1 }
+\sum_{p = 1}^P\left(c^i_p \frac{\partial\phi^{\omega_i}_p}{\partial\nu} -c^{i+1}_p\frac{\partial\phi^{\omega_{i+1}}_p}{\partial\nu}\right)|_{\Gamma_i} = 0
\ {\rm on} \ \Gamma_i \ { \rm for} \ 1< i < I.
\end{split}
\end{equation}

\subsection{\high{The linear} system}
\label{sec:sys}
Once the representation of the solution has been determined and the boundary integral equations \high{are} 
derived, the unknown densities, periodicity constants $c^i_j$ for the proxy surfaces and the 
coefficients of the Rayleigh-Bloch expansion need to be approximated.  This is done by approximating
the boundary integral equations, {\color{black}enforcing} the quasi-periodicity of the solution and 
its flux on the left and right walls, and {\color{black}enforcing} the continuity of the solution 
through the top and bottom of the unit cell.  

In this paper, the boundary integral equations are discretized via a {Nystr\"om} method but the fast 
direct solver can be applied \yzcmt{to} the linear system arising from other discretizations.  Let $N_l$ denote
the number of discretization points on interface $\Gamma_l$.  
As in \cite{Cho:15}, the quasi-periodicity is enforced at points that lie on Gaussian panels between each 
interface on the left and right walls of the unit cell.  Let $M_w$  denote the number of points used 
to enforce periodicity in a layer. (For simplicity
of presentation, we assume this number is the same for all the layers.)
Lastly, the continuity of the integral representation and the Rayleigh-Bloch expansions is enforced at 
collection of $M$ uniformly distributed points on the top and bottom of the unit cell.  The 
Rayleigh-Bloch expansions are truncated at $\pm K$.

The rectangular linear system that arises from these choices has the form

\begin{equation}
 \mthree{\mtx{A}}{\mtx{B}}{\vct{0}}{\mtx{C}}{\mtx{Q}}{\vct{0}}{\mtx{Z}}{\mtx{V}}{\mtx{W}}\vthree{\vct{\hat{\sigma}}}{\vct{c}}{\vct{a}}
 = \vthree{\vct{f}}{\vct{0}}{\vct{0}}
 \label{eq:bigblock}
\end{equation}
where $\mtx{A}$ is a matrix of size $2N\times 2N$ where $N=\sum_{l=1}^I N_l$,
$\mtx{B}$ is a matrix of size $2N\times P$ where $P = \sum_{l=1}^{I+1} P_l$, 
$\mtx{C}$ is a matrix of size \yzcmt{$2(I+1) M_w\times 2N$},
$\mtx{Q}$ is a matrix of size $2(I+1)M_w\times P$,
$\mtx{Z}$ is a matrix of size $4M\times 2N$,
$\mtx{V}$ is a matrix of size $4M\times P$,
and $\mtx{W}$ is a matrix of size $4M\times 2(2K+1)$.
The first row equation enforces the continuity of the scattered field and 
its flux through the interfaces. The second
row equation enforces the quasi-periodicity of the solution and the flux. The last row 
equation enforces continuity
of the integral representation and the Rayleigh-Bloch expansions.  

When the interface geometries 
are complex, a large number of discretization points $N$ are needed to achieve a desired accuracy. 
{\color{black} Because the number of discretization points \yzcmt{on an interface} $N_i$ is
significantly larger \yzcmt{than $M_w$,} $M$, $K$, and $P$ in this scenario,}  
the cost of inverting a matrix the size of $\mtx{A}$ dominates the cost of building 
a direct solver.  For this reason, it is best to build the direct solver in a manner that allows
for the bulk of the computational cost associated with matrices of the size $2N\times 2N$ to be
reused.  \high{We} choose to build a fast 
direct solver for (\ref{eq:bigblock}) via the following block solve:
\begin{align}
 \label{eq:sigma} \vct{\hat{\sigma}} & = -\mtx{A}^{-1}\left([\mtx{B}\ \ \mtx{0}]\vtwo{\vct{c}}{\vct{a}} +\mtx{A}^{-1}\vct{f}\right)\\
 \label{eq:other} \vtwo{\vct{c}}{\vct{a}} & = -\left(\mtwo{\mtx{Q}}{\vct{0}}{\mtx{V}}{\mtx{W}}-\vtwo{\mtx{C}}{\mtx{Z}}\mtx{A}^{-1}[\mtx{B} \ \  \mtx{0}]\right)^\dagger
\vtwo{\mtx{C}}{\mtx{Z}}\mtx{A}^{-1}\vct{f} 
\end{align}
where $\dagger$ denotes the Penrose pseudo-inverse.  

\high{\begin{remark}
 A linear scaling direct solver can be built by processing the block solve in the same order as in \cite{Cho:15}; i.e. solving for $[\vct{c} \ {\vct{a}}]^T$ 
 first.  
 The matrix that needs to be inverted in order to solve for $\vct{\hat{\sigma}}$ is an approximation of 
 the quasi-periodic Green's function and thus can be ill-conditioned when the incident angle is a Wood's anomaly.
\end{remark}}

Each of the matrices \yzcmt{in} (\ref{eq:bigblock}) has a sparsity pattern that can be used to accelerate
the block solve.   The bulk of the acceleration {\color{black}comes} from a fast direct solver for the 
matrix $\mtx{A}$ (see section \ref{sec:fast_Ainv}).

The matrix $\mtx{A}$ is \yzcmt{block tridiagonal}.  
The diagonal blocks of $\mtx{A}$ denoted \yzcmt{by} $\mtx{A}_{ii}$ can be written as the sum 
of two matrices {\color{black}$\mtx{A}^s_{ii}$} and {\color{black}$\mtx{A}^{pm}_{ii}$} where {\color{black}$\mtx{A}^s_{ii}$}
corresponds to
the integral operator on $\Gamma_i$ in the unit cell evaluated on $\Gamma_i$, i.e.
$$ 
\yzcmt{\mtx{A}^s_{ii}=}
\mtwo{-I + D^{\omega_i}_{\Gamma_i,\Gamma_i}-D^{\omega_{i+1}}_{\Gamma_i,\Gamma_i}}
{S^{\omega_{i}}_{\Gamma_i,\Gamma_i}-S^{\omega_{i+1}}_{\Gamma_i,\Gamma_i}}
{T^{\omega_i}_{\Gamma_i,\Gamma_i}-T^{\omega_{i+1}}_{\Gamma_i,\Gamma_i}}{I+D^{*,\omega_i}_{\Gamma_i,\Gamma_i}-D^{*,\omega_{i+1}}_{\Gamma_i,\Gamma_i}},$$
 where $I$ denotes the identity operator, and $\mtx{A}^{pm}_{ii}$ is the contributions from the left and right neighboring 
\yzcmt{copies}, i.e.
$$
\yzcmt{\mtx{A}^{pm}_{ii}=}
\mtwo{D^{\omega_i,pm}_{\Gamma_i,\Gamma_i}-D^{\omega_{i+1},pm}_{\Gamma_i,\Gamma_i}}
{{\color{black}S^{\omega_i,pm}_{\Gamma_i,\{\Gamma_i\}}-S^{\omega_{i+1},pm}_{\Gamma_i,\Gamma_i}}}
{{\color{black}T^{\omega_i,pm}_{\Gamma_i,\Gamma_i}-T^{\omega_{i+1},pm}_{\Gamma_i,\Gamma_i}}}{D^{*,\omega_i,pm}_{\Gamma_i,\Gamma_i}-D^{*,\omega_{i+1},pm}_{\Gamma_i,\Gamma_i}},$$
for $i = 1,\ldots,I$.
The upper diagonal block $\mtx{A}_{i,i+1}$ corresponds to the integral operators on $\Gamma_{i+1}$ being 
evaluated on $\Gamma_i$, i.e.
$$
\yzcmt{\mtx{A}_{i,i+1}=}
\mtwo{-\tilde{D}^{\omega_{i+1}}_{\Gamma_i,\Gamma_{i+1}}}
{-\tilde{S}^{\omega_{i+1}}_{\Gamma_i,\Gamma_{i+1}}}
{-\tilde{T}^{\omega_{i+1}}_{\Gamma_i,\Gamma_{i+1}}}
{-\tilde{D}^{*,\omega_{i+1}}_{\Gamma_i,\Gamma_{i+1}}},$$
for $i = 1,\ldots,I-1$.
 The lower diagonal blocks $\mtx{A}_{i,i-1}$ correspond\high{s} to the integral operators on $\Gamma_{i-1}$ 
being evaluated on $\Gamma_i$, i.e.
$$
\yzcmt{\mtx{A}_{i,i-1}=}
\mtwo{\tilde{D}^{\omega_{i}}_{\Gamma_i,\Gamma_{i-1}}}
{\tilde{S}^{\omega_{i}}_{\Gamma_i,\Gamma_{i-1}}}
{\tilde{T}^{\omega_{i}}_{\Gamma_i,\Gamma_{i-1}}}
{\tilde{D}^{*,\omega_{i}}_{\Gamma_i,\Gamma_{i-1}}},$$
for $i= 2,\ldots,I$.

The matrix $\mtx{B}$ is upper block diagonal with block \yzcmt{defined} by 
\begin{equation} 
\mtx{B}_{i,i} = \left[ \begin{array}{ccc} \phi_1^{\omega_i}|_{\Gamma_i} & \cdots & \phi^{\omega_i}_P|_{\Gamma_i} \\
                          \frac{\partial\phi_1^{\omega_i}}{\partial \vct{n}}|_{\Gamma_i} & \cdots & \frac{\partial \phi^{\omega_i}}{\partial \vct{n}}_P|_{\Gamma_i}
                         \end{array}\right] \ {\rm and} \
\mtx{B}_{i,i+1} = \left[ \begin{array}{ccc} -\phi_1^{\omega_{i+1}}|_{\Gamma_i} & \cdots & -\phi^{\omega_{i+1}}_P|_{\Gamma_i} \\
                          -\frac{\partial\phi_1^{\omega_{i+1}}}{\partial \vct{n}}|_{\Gamma_i} & \cdots & -\frac{\partial \phi^{\omega_{i+1}}}{\partial \vct{n}}_P|_{\Gamma_i}
                         \end{array}\right]                         
                         \label{eq:Bmat}
                         \end{equation}

for $i = 1,\ldots,I$.
The matrix $\mtx{C}$ is lower block diagonal with blocks \yzcmt{defined} by
\begin{equation}
\mtx{C}_{i,i} = \mtwo{\alpha^{-2}D^{\omega_i}_{R_i+\vct{d},\Gamma_i}-\alpha D^{\omega_i}_{L_i-\vct{d},\Gamma_i}}
                  {\alpha^{-2}S^{\omega_i}_{R_i+\vct{d},\Gamma_i}-\alpha S^{\omega_i}_{L_i-\vct{d},\Gamma_i}}
                  {\alpha^{-2}T^{\omega_i}_{R_i+\vct{d},\Gamma_i}-\alpha T^{\omega_i}_{L_i-\vct{d},\Gamma_i}}
                  {\alpha^{-2}D^{*,\omega_i}_{R_i+\vct{d},\Gamma_i}-\alpha D^{*,\omega_i}_{L_i-\vct{d},\Gamma_i}} 
                  \ {\rm and}                  \label{eq:Cmat1}
\end{equation}

                   \begin{equation}
                  \mtx{C}_{i,i-1} = \mtwo{\alpha^{-2}D^{\omega_i}_{R_i+\vct{d},\Gamma_{i-1}}-\alpha D^{\omega_i}_{L_i-\vct{d},\Gamma_{i-1}}}
                  {\alpha^{-2}S^{\omega_i}_{R_i+\vct{d},\Gamma_{i-1}}-\alpha S^{\omega_i}_{L_i-\vct{d},\Gamma_{i-1}}}
                  {\alpha^{-2}T^{\omega_i}_{R_i+\vct{d},\Gamma_{i-1}}-\alpha T^{\omega_i}_{L_i-\vct{d},\Gamma_{i-1}}}
                  {\alpha^{-2}D^{*,\omega_i}_{R_i+\vct{d},\Gamma_{i-1}}-\alpha D^{*,\omega_i}_{L_i-\vct{d},\Gamma_{i-1}}}  
                  \label{eq:Cmat2}
\end{equation}

for $i = 1,\ldots, I$ and $i = 2,\ldots,I+1$, respectively.                  
The matrix $\mtx{Q}$ is block diagonal with blocks given by 
\begin{equation}
\mtx{Q}_{ii} = \left[\begin{array}{ccc}
           \alpha^{-1}\phi_1^{\omega_i}|_{R_i} -\phi_1^{\omega_i}|_{L_i} & \cdots &   \alpha^{-1}\phi_P^{\omega_i}|_{R_i} -\phi_P^{\omega_i}|_{L_i}\\
           \alpha^{-1}\frac{\partial\phi_1^{\omega_i}}{\partial \vct{n}}|_{R_i} -\frac{\partial\phi_1^{\omega_i}}{\partial \vct{n}}|_{L_i} & \cdots & 
           \alpha^{-1}\frac{\partial\phi_P^{\omega_i}}{\partial \vct{n}}|_{R_i} -\frac{\partial\phi_P^{\omega_i}}{\partial \vct{n}}|_{L_i}\\
                       \end{array}\right] 
                       \label{eq:Qmat}
\end{equation}

                  for $i = 1,\ldots,I+1$.  

The matrices $\mtx{Z}$, $\mtx{V}$, and $\mtx{W}$ are sparse matrices of the form
$$\mtx{Z}= \left[\begin{array}{cccc}
         \mtx{Z}_U & \mtx{0} & \cdots & \mtx{0}\\
         \mtx{0} & \cdots &\mtx{0} &\mtx{Z}_D
        \end{array}\right], \ \mtx{V} = \left[\begin{array}{cccc}
         \mtx{V}_U & \mtx{0} & \cdots & \mtx{0}\\
         \mtx{0} & \cdots &\mtx{0} &\mtx{V}_D
        \end{array}\right], \ {\rm and} \ \mtx{W} = \mtwo{\mtx{W}_U}{\mtx{0}}{\mtx{0}}{\mtx{W}_D}$$
where 
\begin{equation}\mtx{Z}_U = \mtwo{\tilde{D}^{\omega_1}_{U,\Gamma_1}}{\tilde{S}^{\omega_1}_{U,\Gamma_1}}{\tilde{T}^{\omega_1}_{U,\Gamma_1}}{\tilde{D}^{*,\omega_1}_{U,\Gamma_1}}, \
\mtx{Z}_D = \mtwo{\tilde{D}^{\omega_{I+1}}_{D,\Gamma_I}}{\tilde{S}^{\omega_{I+1}}_{D,\Gamma_I}}{\tilde{T}^{\omega_{I+1}}_{U,\Gamma_I}}
{\tilde{D}^{*,\omega_{I+1}}_{U,\Gamma_I}},
\label{eq:Zmat}
\end{equation}
\begin{equation}\mtx{V}_U = \left[\begin{array}{ccc}{\phi^{\omega_1}_1|_U}& \cdots &\phi^{\omega_1}_P|_{U}\\
                   {\frac{\phi^{\omega_1}_1}{\partial \nu}|_U}& \cdots &\frac{\phi^{\omega_1}_P}{\partial\nu}|_{U}
                  \end{array}\right], \
                  \mtx{V}_D = \left[\begin{array}{ccc}{\phi^{\omega_{I+1}}_1|_D}& \cdots &\phi^{\omega_{I+1}}_P|_{D}\\
                   {\frac{\phi^{\omega_{I+1}}_1}{\partial \nu}|_D}& \cdots &\frac{\phi^{\omega_{I+1}}_P}{\partial\nu}|_{U}
                  \end{array}\right], \label{eq:Vmat}\end{equation}
                  
                  \begin{equation}\mtx{W}_{U} = \left[\begin{array}{ccc}
                                         -e^{i\kappa_{-K}x}|_U & \cdots & -e^{i\kappa_{K}x}|_U \\
                                         -ik^U_{-K}e^{i\kappa_{-K}x}|_U&\cdots & -ik^U_{K}e^{i\kappa_{K}x}|_U 
                                        \end{array}\right], \ {\rm and} \                                         
\mtx{W}_{D} = \left[\begin{array}{ccc}
                                         -e^{i\kappa_{-K}x}|_D & \cdots & -e^{i\kappa_{K}x}|_D \\
                                         ik^D_{-K}e^{i\kappa_{-K}x}|_D&\cdots & ik^D_{K}e^{i\kappa_{K}x}|_D 
                                        \end{array}\right].
\label{eq:Wmat}\end{equation}

The matrices $\mtx{W}_U$ and $\mtx{W}_D$ correspond to the evaluation of the terms in the 
Rayleigh-Bloch expansions \yzcmt{at points on the top and bottom of the unit cell 
where continuity of the solution is enforced.}

\section{\high{The fast} direct solver}
\label{sec:FDS}
While exploiting the sparsity of the matrices can accelerate the construction
of a direct solver, the speed gains are not sufficient for applications when 
the interface geometries are complex.  {\color{black}When the interface
geometries are complex, the cost of building a direct solver for the rectangular
system is dominated by the cost of inverting $\mtx{A}$.}  The fast direct solver proposed in 
this section exploits not only \yzcmt{the} sparsity but also the data sparse nature of \yzcmt{the matrix 
$\mtx{A}$}.  

The \yzcmt{foundation} of the fast direct solver is a fast inversion technique
for $\mtx{A}$ presented in section \ref{sec:fast_Ainv}. 
The fast inversion of $\mtx{A}$ allows for $\vct{\hat{\sigma}}$ to be computed for a cost that 
scales linearly \yzcmt{with respect to $N$} via equation (\ref{eq:sigma}).  Constructing and applying {\color{black}
an approximation of }the pseudo-inverse of the Schur complement 
 \begin{equation}
 \yzcmt{
\mtx{S}= -\left(\mtwo{\mtx{Q}}{\vct{0}}{\mtx{V}}{\mtx{W}}-\vtwo{\mtx{C}}{\mtx{Z}}\mtx{A}^{-1}[\mtx{B} \ \  \mtx{0}]\right)}
\label{eq:Schur}  
 \end{equation}
is needed to find $\vct{c}$ and $\vct{a}$ via (\ref{eq:other}).  
{\color{black}The approximate pseudo-inverse is created by first computing an 
$\epsilon_{\rm Schur}$\textcolor{black}{-truncated} singular value decomposition (SVD) 
and \yzcmt{then} applying the pseudo-inverse of this factorization.}

\begin{definition}
Let $\mtx{U}\mtx{\Sigma}\mtx{T}^*$ be the SVD of \yzcmt{the} Schur complement matrix $\mtx{S}$ of size $(2(I+1)M_w+4M)\times (P+2(2K+1))$ where 
$\mtx{\Sigma}$ is a diagonal rectangular matrix with entries of the singular values of $\mtx{S}$ and
matrices $\mtx{U}$ and $\mtx{T}$ are unitary matrices of size $(2(I+1)M_w+4M)\times(2(I+1)M_w+4M)$ and $(P+2(2K+1))\times(P+2(2K+1))$, respectively.  
Then the $\epsilon_{\rm Schur}$
\textcolor{black}{-truncated} SVD
 is 
$$\mtx{\hat{U}}\mtx{\hat{\Sigma}} \mtx{\hat{T}}^*$$
where $\mtx{\hat{\Sigma}}$ is a diagonal square matrix of size $l\times l$ where $l$ is the
number of singular values of $\mtx{S}$ that are larger than $\epsilon_{\rm Schur}$, $\mtx{\hat{U}}$ is 
an $(2(I+1)M_w+4M)\times l$ matrix and $\mtx{\hat{T}}$ is an $(P+2(2K+1))\times l$ matrix.
\end{definition}
\textcolor{black}{In practice, we found $\epsilon_{Schur}=10^{-13}$ is a good choice when the desired accuracy for the solution is $10^{-10}$.}

Then $\vct{c}$ and $\vct{a}$ can be approximated by 
$$\vtwo{\vct{c}}{\vct{a}}\approx \mtx{\hat{T}} \mtx{\hat{\Sigma}}^{-1} \mtx{\hat{U}}^* \vtwo{\mtx{C}}{\mtx{Z}}\mtx{A}^{-1}\vct{f}.$$

The most efficient way to find $\vct{c}$ and $\vct{a}$ is to 
apply the matrices from right to left in this equation meaning that the 
vectors are found via a collection \yzcmt{of} matrix vector multiplies. 

\begin{remark}
The cost of constructing the truncated SVD for $\mtx{S}$ scales cubically with respect to the number of interfaces $I$ but 
is constant with respect to the number of points on the interfaces. 
\end{remark}

\high{
Combining the fact that many of the matrices (less scalar factors) in (\ref{eq:bigblock}) can be re-used
for multiple incident angles (see section \ref{sec:reuse})
with a fast direct solver for $\mtx{A}$} results 
in a fast direct solver that is ideal for problems where many solves are required.  An
additional key feature of the fast direct
solver is that the bulk of the precomputation can be re-used if an interface $\Gamma_j$ or a
wave number $\omega_j$ is changed (see Section \ref{sec:extensions}).

\subsection{Fast inversion of $\mtx{A}$}
\label{sec:fast_Ainv} 
The key to building the fast direct solver for the block system (\ref{eq:bigblock}) is 
having a fast way of inverting $\mtx{A}$.  This technique is designed to make solves 
for different Bloch phases as efficient as possible.  

The solver considers the matrix $\mtx{A}$ written as the sum of two matrices
\begin{equation}
 \label{eq:tri_sum}
 \begin{split}
 \mtx{A} &= \underbrace{\left[ \begin{array}{ccccc}
                   {\color{black}\mtx{A}^s_{11}} &\mtx{0}   &\mtx{0}&\mtx{0}&\mtx{0}\\
                   \mtx{0} & {\color{black}\mtx{A}^s_{22}}  &\mtx{0}&\mtx{0} &\mtx{0}\\
                   \mtx{0}      & \mtx{0} & \ddots   & \mtx{0} & \mtx{0} \\
                   \mtx{0}      & \mtx{0}   & \mtx{0}&{\color{black}\mtx{A}^s_{ (N-1)(N-1)}} &\mtx{0}\\
                   \mtx{0}      & \mtx{0}   & \mtx{0} &\mtx{0} &{\color{black}\mtx{A}^s_{NN}}\\
                  \end{array}\right]}_{\mtx{A}_0}\\
                  &+\underbrace{
                  \left[ \begin{array}{ccccc}
                   {\color{black}\mtx{A}^{pm}_{11}} & \mtx{A}_{12} &\mtx{0}&\mtx{0}&\mtx{0}\\
                   \mtx{A}_{21} & {\color{black}\mtx{A}^{pm}_{22}} & \mtx{A}_{23} &\mtx{0} &\mtx{0}\\
                   \mtx{0}      & \mtx{0} & \ddots   & \ddots & \mtx{0} \\
                   \mtx{0}      & \mtx{0}   &\mtx{A}_{(N-1),(N-2)} &{\color{black}\mtx{A}^{pm}_{(N-1)(N-1)}} &\mtx{A}_{(N-1),N}\\
                   \mtx{0}      & \mtx{0}   & \mtx{0} &\mtx{A}_{N,(N-1)} &{\color{black}\mtx{A}^{pm}_{NN}}\\
                  \end{array}\right]}_{\hat{\mtx{A}}}
                  \end{split}
\end{equation}
where the block diagonal matrix $\mtx{A}_0$ {\color{black}whose entries are self-interaction matrices}
and $\mtx{\hat{A}}$ is the block tridiagonal matrix where the diagonal blocks correspond to the 
interaction of an interface with its left and right neighbors and the off-diagonal blocks correspond
to the interactions between the interfaces directly above and below each other.  
  Since the submatrices in $\hat{\mtx{A}}$ correspond to ``far'' interactions,
they are numerically low rank. Let $\mtx{L}\mtx{R}$ denote the low rank factorization of $\hat{\mtx{A}}$
where $\mtx{L}$ and $\mtx{R}^T$ are $2N\times k_{\rm tot}$ matrices and $k_{\rm tot}$ is the numerical 
rank of $\hat{\mtx{A}}$.
Section \ref{sec:lowrank} presents a technique for constructing this factorization.  Then $\mtx{A}$ can 
be approximated by 
{\color{black}$$\mtx{A} \approx \mtx{A}_{0} + \mtx{L}\mtx{R}.$$}
The advantage of this representation is that the factors $\mtx{L}$ and $\mtx{R}$ can be computed in 
a way that is independent of Bloch phase as presented in section \ref{sec:lowrank}.  Additionally, 
the inverse can be formulated via a Woodbury formula \cite{golub}

{\color{black}\begin{equation}\mtx{A}^{-1} \approx (\mtx{A}_{0} +\mtx{L}\mtx{R} )^{-1} = \mtx{A}_{0}^{-1} -
\mtx{A}_{0}^{-1}\mtx{L}\left(\mtx{I}+\mtx{R}\mtx{A}_{0}^{-1}\mtx{L}\right)^{-1}\mtx{R}\mtx{A}_{0}^{-1}.\label{eq:inv}\end{equation}}

Not only is the matrix $\mtx{A}_0$ block diagonal but each of the diagonal blocks is amenable to 
a fast direct solver such as \textit{Hierarchically
Block Separable (HBS)} methods \cite{2012_martinsson_FDS_survey,2012_ho_greengard_fastdirect,2013_3DBIE}
which are closely related to the \textit{Hierarchically Semi-Separable (HSS)}
 \cite{2009_xia_superfast,2007_shiv_sheng,2010_xia}, the Hierarchical interpolative
 factorization \cite{2014_HIF}, the $\mathcal{H}$ and $\mathcal{H}^{2}$-matrix methods
 \cite{2010_borm_book,2004_borm_hackbusch}.  
 \high{Thus an approximate inverse of $\mtx{A}_0$ can be constructed and applied} for a
 cost that scales linearly \yzcmt{with respect to} the number of discretization points on the interfaces.
 This computation is independent of Bloch phase.  
 \textcolor{black}{The condition number of $(\mtx{I}+\mtx{R}\mtx{A}_0^{-1}\mtx{L})$ is bounded above by the 
 product of the condition number of $\mtx{A}_0$ and $(\mtx{A}_0+\mtx{L}\mtx{R})$ \cite{Yip_woodburystability}.
Since both $\mtx{A}$ and $\mtx{A}_0$ result from the discretization of 
second kind boundary integral equations, they are well-conditioned.  Thus 
 applying the Woodbury formula is numerically stable.
}

It is never necessary to construct the approximation of $\mtx{A}^{-1}$.  It is 
only necessary to have a fast algorithm for applying it to a vector $\vct{f}\in \mathbb{C}^{2N\times 1}$,
i.e. a fast algorithm is needed for evaluating 
{\color{black}\begin{equation}\mtx{A}^{-1}\vct{f}\approx  \mtx{A}_{0}^{-1}\vct{f} -
\mtx{A}_{0}^{-1}\mtx{L}\left(\mtx{I}+\mtx{R}\mtx{A}_{0}^{-1}\mtx{L}\right)^{-1}\mtx{R}\mtx{A}_{0}^{-1}\vct{f}\label{eq:invapply}.\end{equation}}
The fast direct solver for $\mtx{A}_0$ and the block structure of the matrices $\mtx{L}$ and $\mtx{R}$
allow for $\mtx{A}_{0}^{-1}\mtx{L}$ and 
$\mtx{R}\mtx{A}_{0}^{-1}\vct{f}$ to be evaluated for a cost that scales linearly with $N$.
Thanks to the sparsity pattern of the matrices, the intermediate matrix
$\mtx{S}_2 = \mtx{I}+\mtx{R}\mtx{A}_{0}^{-1}\mtx{L}$ 
of size $k_{\rm tot}\times k_{\rm tot}$ 
that needs to be inverted is block tridiagonal.  
 Appendix \ref{app:1} reports 
on the construction of $\mtx{S}_2$.  
Since $k_{\rm tot}$ is much smaller than $N$ in practice, the inverse of $\mtx{S}_2$ can be applied rapidly 
using a block variant of the Thomas algorithm.
This computation needs to be done for each new Bloch phase since $\mtx{L}$ and $\mtx{R}$ are 
dependent on Bloch phase. 

\begin{remark}
To achieve nearly optimal ranks in the construction of the fast direct solver, it 
is advantageous to reorder the matrices in $\mtx{A}$ according to the
 physical location of the unknowns.  For example, if there are $N_1$ discretization 
 points on $\Gamma_1$, the unknowns are $\sigma_{1,1}, \ldots, \sigma_{1,N_1}$ and
 $\tau_{1,1},\ldots, \tau_{1,N_1}$, etc. Then the matrices should be 
 ordered so $\vct{\hat{\sigma}}$  is as follows
 $$\vct{\hat{\sigma}}^T = \left[ \sigma_{1,1}, \tau_{1,1}, \cdots, \sigma_{1,N_1}, \tau_{1,N_1},
 \cdots, \sigma_{I,1}, \tau_{I,1}, \cdots, \sigma_{I,N_I}, \tau_{I,N_I}\right].$$
\end{remark}

\subsubsection{Low rank factorization of $\hat{\mtx{A}}$}
\label{sec:lowrank}

The technique for creating the low rank factorizations of the blocks in $\hat{\mtx{A}}$ 
is slightly different depending on whether or not the block is {\color{black}on the} diagonal.  
This section begins by presenting the technique for creating low rank factorizations of the 
diagonal blocks.  Then the technique for creating the low rank \yzcmt{factorizations} of 
the off-diagonal blocks is presented.  

Recall the diagonal blocks of $\hat{\mtx{A}}$ are {\color{black}$\mtx{A}^{pm}_{ii}$} 
\yzcmt{corresponding} to the discretized
version of 
$$
\displaystyle({\color{black}\mathcal{S}^{\omega,pm}_{\Gamma_i,\Gamma_i}}\rho)(\vct{x}) = 
\alpha \underbrace{ \int_{\Gamma_i} G_\omega(\vct{x},\vct{y}+\vct{d})\rho(\vct{y})dl(\vct{y})}_{
\mbox{right copy: } \displaystyle({\color{black}\mathcal{S}^{\omega,p}_{\Gamma_i,\Gamma_i}}\rho)(\vct{x})}
+
\alpha^{-1} \underbrace{ \int_{\Gamma_i} G_\omega(\vct{x},\vct{y}-\vct{d})\rho(\vct{y})dl(\vct{y})}_{
\mbox{left copy: } \displaystyle({\color{black}\mathcal{S}^{\omega,m}_{\Gamma_i,\Gamma_i}}\rho)(\vct{x})}.
$$
{\color{black} The matrix $\mtx{A}^{pm}_{ii}$ can be written as the sum of two matrices that are independent of 
Bloch phase; $\mtx{A}^{pm}_{ii} = \alpha\mtx{A}^p_{ii}+\alpha^{-1} \mtx{A}^m_{ii}.$}
Thus by creating low rank factorizations of {\color{black}$\mtx{A}^p_{ii}$} and {\color{black}$\mtx{A}^m_{ii}$} independently,
the factorizations can be used for any Bloch phase $\alpha$.  Let {\color{black} $\mtx{L}^p_{i} \mtx{R}^p_{i}$ }
and {\color{black}$\mtx{L}^m_{i} \mtx{R}^m_{i} $} denote the low rank approximations of {\color{black}$\mtx{A}^p_{ii}$}
and {\color{black}$\mtx{A}^m_{ii}$} respectively.  These two approximations are combined to create 
a low rank approximation of {\color{black}$\mtx{A}^{pm}_{ii}$} as follows:
$${\color{black}\mtx{A}^{pm}_{ii}\approx
{\color{black}\underbrace{
[ \mtx{L}^p_{i},\; \mtx{L}^m_{i}]}_{\mtx{L}^{pm}_{ii}}}}
\underbrace{
\begin{bmatrix}
\alpha{\color{black}\mtx{R}^p_{i}} \\
\alpha^{-1}{\color{black}\mtx{R}^m_{i}} 
\end{bmatrix}}_{{\color{black}\mtx{R}^{pm}_{ii}}}
$$

{\color{black} The technique used to create the low rank factorizations is similar to
the one used in \yzcmt{\cite{2016_periodic_stokes}.  The new technique}
has an extra step to keep the rank $k_{\rm tot}$ small.}

For brevity,this manuscript \yzcmt{only presents the technique for compressing the interaction 
with the left neighbor (i.e. computing the low rank factorization of 
{\color{black}$\mtx{A}^m_{ii}$}).}  
The technique for compressing the 
interaction with the right neighbor follows directly.

We choose to build the factorization via the interpolatory decomposition \cite{gu1996,lowrank}
defined as follows.

\begin{definition}
 The \textit{interpolatory decomposition} of a $m\times n$ matrix $\mtx{M}$ that has rank $l$ is
the factorization 
$$ \mtx{M} = \mtx{P}\mtx{M}(J(1:l),:)$$
where $J$ is a vector of integers $j$ such $1\leq j\leq m$, and $\mtx{P}$ is a $m\times l$ matrix that contains a $l \times l$ identity matrix.
Namely, $\mtx{P}(J(1:l),:) = \mtx{I}_l$.  
\end{definition}
%

Creating the low rank factorization of {\color{black}$\mtx{A}^m_{ii}$} by directly plugging it into the interpolatory 
decomposition has a computational cost of $O(N_i^2k_i)$ where $k_i$ is the 
numerical rank of {\color{black}$\mtx{A}^m_{ii}$}.  This would result in a solution
technique that has a computational cost that scales quadratically, \emph{not linearly}, with respect
to $N$.  To {\color{black}achieve} the linear computational complexity,
we utilize potential theory.

\begin{figure}[h]
\centering
\begin{tabular}{c}
{\begin{picture}(200,150)
\put(-30,-10){\includegraphics[trim={15cm 2cm 4cm 2cm},clip,scale=0.5]{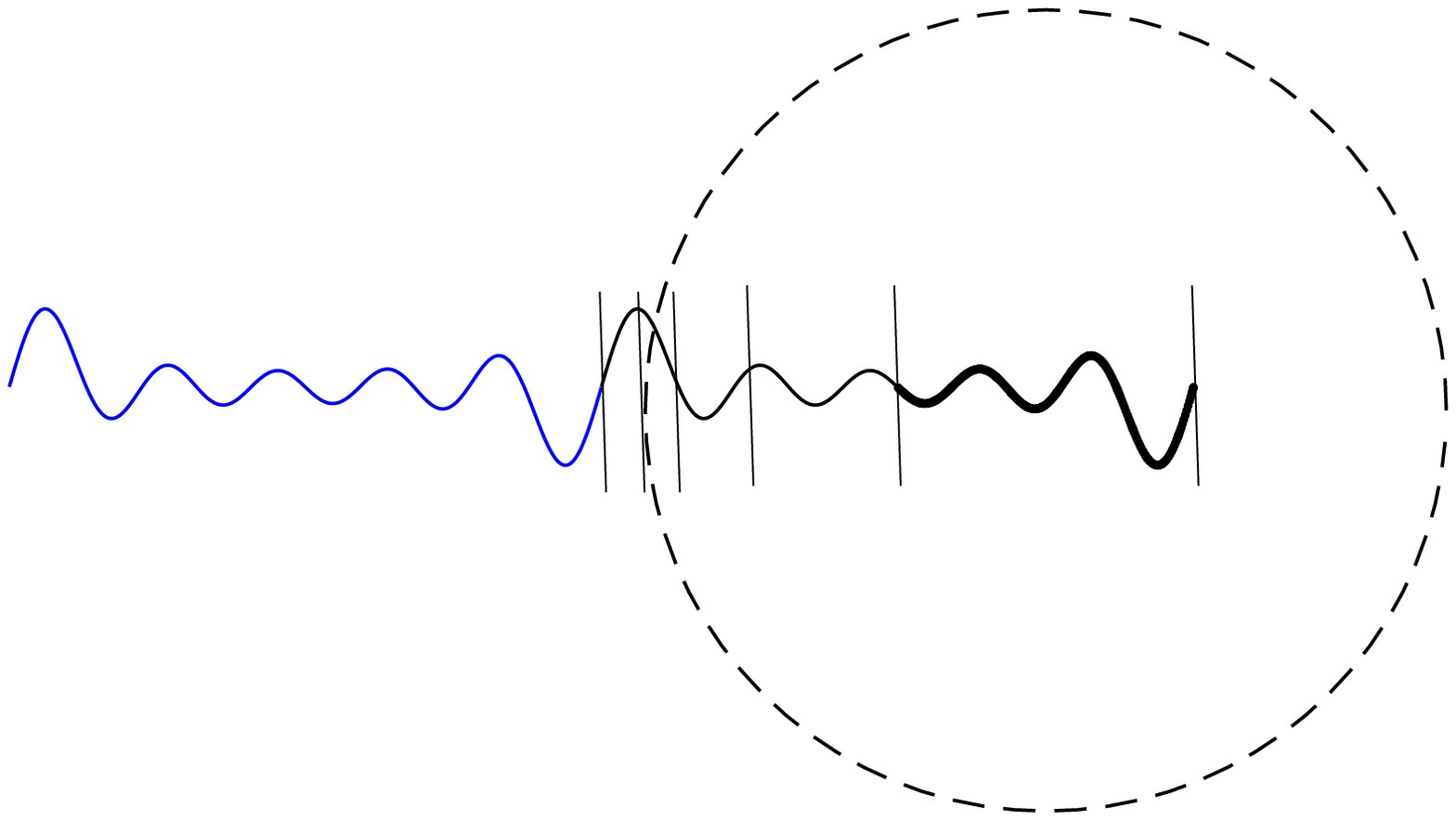}}
\put(40,60){\textcolor{black}{$\Gamma^m_{i}$}}
\put(180,65){$\gamma_1$}
 \end{picture}
}\\
(a)\\
{\begin{picture}(200,150)
\put(10,35){  \includegraphics[trim={18cm 4cm 15cm 5cm},clip,scale=0.4]{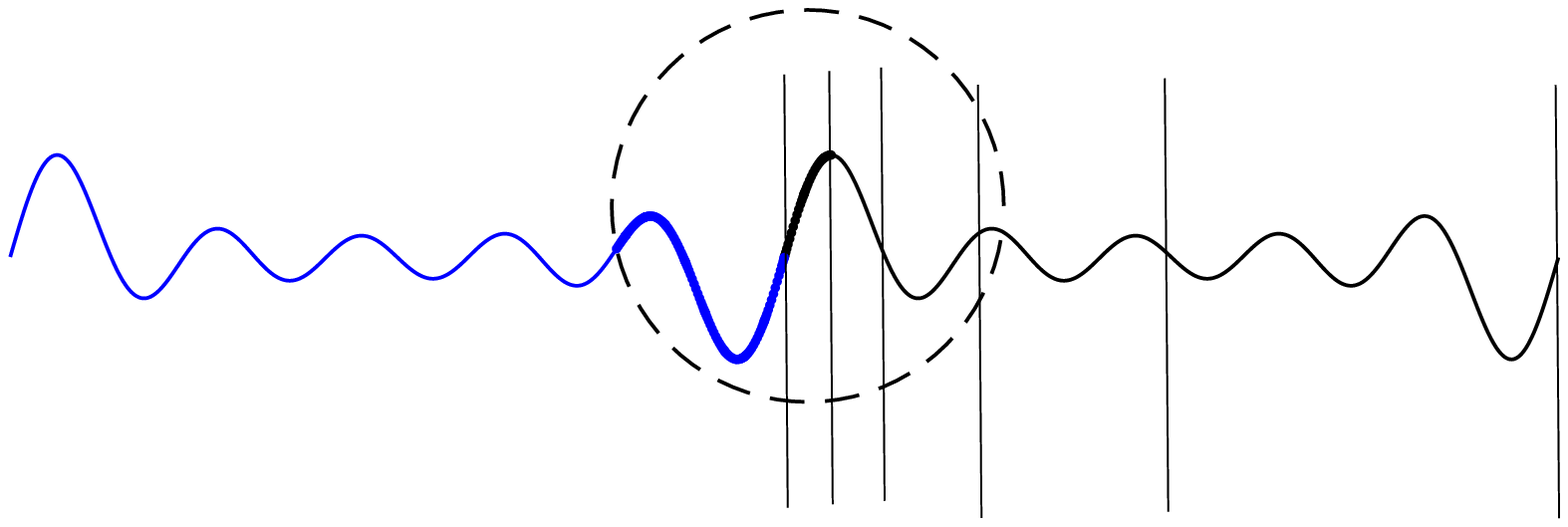}}
\put(30,80){\textcolor{black}{$\Gamma^m_{i}$}}
\put(115,90){$\gamma_5$}
 \end{picture}
}\\
(b)
\end{tabular}
 \caption{Illustration of the dyadic refinement partitioning of \yzcmt{$\Gamma_i$} with 5 levels
 refinement and geometries for compressing {\color{black}$\mtx{A}^m_{ii}$}. 
 (a) Illustration of the proxy surface (dashed circle) used to compress neighbor interactions when $\gamma_l$ is far
 \yzcmt{from $\Gamma_i^m$}. 
 (c) Illustration of the proxy surface (dashed circle) and near points (bold \yzcmt{blue} curve on {\color{black}$\Gamma^m_{i}$}) when 
 $\gamma_l$ is touching {\color{black}$\Gamma^m_{i}$}. }
 \label{fig:neigh} 
\end{figure}

Recall $\Gamma_i$ denotes the part of the $i^{\rm th}$ interface in the unit cell.  Let {\color{black}$\Gamma^m_{i}$}
denote the part of the $i^{\rm th}$ interface in the left neighboring cell.  First $\Gamma_i$ is partitioned
into a collection of $S$ segments $\gamma_j$ via dyadic refinement where the segments get smaller as they approach
{\color{black}$\Gamma^m_{i}$} so that $\Gamma_i = \displaystyle\cup_{j=1}^S \gamma_j$. Figure \ref{fig:neigh} illustrates a partitioning
when compressing the interaction of $\Gamma_i$ with {\color{black}$\Gamma^m_{i}$}.   The refinement is stopped
when the segment closest to {\color{black}$\Gamma^m_{i}$} has less than $n_{\rm max}$ points on it.  Typically, $n_{\rm max} = 45$ 
is a good choice.

For each segment $\gamma_j$ not touching {\color{black}$\Gamma^m_{i}$}, consider a circle concentric with the bounding box containing
$\gamma_j$ with a radius slightly less than the distance from the center of the bounding box \yzcmt{to} {\color{black}$\Gamma^m_{i}$}. 
Figure \ref{fig:neigh}(a) illustrates the proxy surface for $\gamma_1$ when there \yzcmt{are} 5 levels of dyadic refinement toward
{\color{black}$\Gamma^m_{i}$}.
From potential theory, we know that any field generated by sources outside of this circle can be approximated
to high accuracy by placing enough equivalent charges on the circle.  In practice, it is enough to place a 
small number of \emph{proxy points} evenly on the circle.  Let $n_{\rm proxy}$ denote the number of 
proxy points on the circle.  
For problems where the direct solver scales linearly, $n_{\rm proxy}$ is 
small and chosen to be a constant independent of $\omega_i$.
   For the experiments in this paper, it is sufficient to {\color{black}set} $n_{\rm proxy} = 80$.
  Let $n_j$ 
denote the number of points on $\gamma_j$.  An interpolatory decomposition can be constructed for the 
matrix $\mtx{A}^{\rm proxy}$ capturing the interaction between $\gamma_j$ and the proxy points.    The 
result is an index vector $J_j$ and an interpolation matrix $\mtx{P}_j$ of size $n_j\times k_j$ where $k_j$ is the 
numerical rank of $\mtx{A}^{\rm proxy}$.  For $\gamma_S$ (the segment 
touching {\color{black}$\Gamma^m_{i}$}), $n_{\rm proxy}$ proxy points are placed uniformly on a circle of radius $1.75$
\yzcmt{times larger than the radius of} 
the smallest circle containing all the points on $\gamma_S$.  
All the points on {\color{black}$\Gamma^m_{i}$} inside the circle are 
labeled \emph{near} points and indexed $I_{\rm near}$.  Figure \ref{fig:neigh}(b) illustrates
the proxy circle and near points for $\gamma_5$ when there \yzcmt{are} 5 levels of dyadic refine toward
{\color{black}$\Gamma^m_{i}$}.  
An interpolatory decomposition is then 
performed on $\left[  {\color{black}\mtx{A}^m_{ii}}(\gamma_S,I_{\rm near}) \ | \  \mtx{A}^{\rm proxy} \right]$.  The result
is an index vector $J_S$ and an interpolation matrix $\mtx{P}_S$ of size $n_S\times k_S$.

The low rank factorization of the matrix {\color{black}$\mtx{A}^m_{ii}$} can be constructed with the result of this compression 
procedure.  Let $J = [J_1(1:k_1),\ldots,J_S(1:k_S)]$
denote an index vector
\yzcmt{consisting of the index vectors for each segment}.
  Then {\color{black}$\mtx{L}^m_{i}$} is a block diagonal matrix with block entries $\mtx{P}_j$ for $J =1,\ldots, S$ 
and {\color{black}$\mtx{R}^m_{i} = \mtx{A}^m_{ii}(J,:)$}.  
\yzcmt{The points on 
$\Gamma_i$ corresponding to the index vector $J$ are called the \emph{skeleton points}.}    

{\color{black}The rank of this factorization is far from optimal} and will result in {\color{black} an excessively large}
constant prefactor in the application of the Woodbury formula (\ref{eq:inv}).  
  {\color{black}A recompression step is necessary to resolve this problem.}  Let $k_{\rm orig}$ denote the 
  rank of the original approximate factorization, i.e. the length of $J$.
If $k_{\rm orig}$ is small enough, applying the interpolatory decomposition to 
{\color{black}$\mtx{A}^m_{ii}(J,:)$} can be done \yzcmt{efficiently resulting in an} index vector $J_{\rm up}$ and interpolation 
matrix $\mtx{P}_{\rm up}$ of size $k_{\rm orig}\times k_{\rm up}$. Let $\mtx{L}_{\rm up}=  \mtx{P}_{\rm up}$.
Otherwise, the interpolatory decomposition can 
be applied to the submatrices corresponding to a lump of the segments at \yzcmt{a} time.  For example, 
suppose $S$ is even, then the segments can bunched two at \yzcmt{a} time.  The interpolatory 
decomposition can be applied to ${\color{black}\mtx{A}^m_{ii}}([J_{j}(1:k_j),J_{j+1}(1:k_{j+1})],:)$ 
for \yzcmt{odd values of $j$}.             
The resulting interpolation matrices are the block entries
for the block diagonal matrix $\mtx{L}_{\rm up}$.  The corresponding index vector $J_{\rm up}$ 
is formed in a similar manner to the vector $J$.  Finally the low rank factorization of 
{\color{black}$\mtx{A}^m_{ii}$} can be formed by multiplying {\color{black}$\mtx{L}^m_{i}$}
\yzcmt{by} $\mtx{L}_{\rm up}$ and using the updated skeleton of $J(J_{\rm up})$.  In other words,
{\color{black}$\mtx{L}^m_{i}= \mtx{L}^m_{i}\mtx{L}_{\rm up}$} and {\color{black}$\mtx{R}^m_{i} = \mtx{A}^m_{ii}(J(J_{\rm up}),:)$}.

The technique for constructing the low rank factorization of the off-diagonal blocks of $\mtx{\hat{A}}$ is 
similar.  Recall that each off-diagonal block  $\mtx{A}_{ij}$, for $i\neq j$, corresponds to the 
discretization of the following integral operator where $\vct{x}\in\Gamma_i$:
\begin{align*}
\yzcmt{
\displaystyle(\tilde{\mathcal{S}}_{\Gamma_i,\Gamma_j}^\omega\rho)(\vct{x})} 
= &
 \int_{\Gamma_j} G_\omega(\vct{x},\vct{y})\rho(\vct{y})dl(\vct{y})
+
\alpha  \int_{\Gamma_j} G_\omega(\vct{x},\vct{y}+\vct{d})\rho(\vct{y})dl(\vct{y})\\
&+
\alpha^{-1} \int_{\Gamma_j} G_\omega(\vct{x},\vct{y}-\vct{d})\rho(\vct{y})dl(\vct{y}).
\end{align*}
It is natural to write $\mtx{A}_{ij}$ as the summation of three parts,
$$
\mtx{A}_{ij}= \mtx{A}_{0,ij}+\alpha {\color{black}\mtx{A}^p_{ij}} +\alpha^{-1} {\color{black}\mtx{A}^m_{ij}},
$$
where $\mtx{A}_{0,ij}$, {\color{black}$\mtx{A}^p_{ij}$}, and {\color{black}$\mtx{A}^m_{ij}$} are the discrete approximations of
the corresponding integral operators.

\begin{figure}[h]
\centering
\begin{picture}(200,200)(0,0)
\put(-45,0){\includegraphics[trim={1cm 1cm 1cm 1cm},clip,scale=0.5]{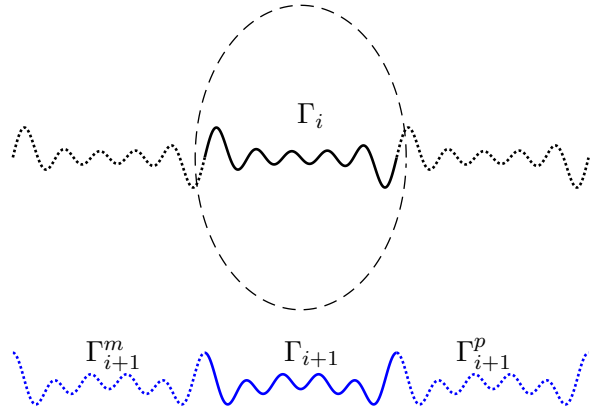}}
\put(85,125){$\Gamma_i$}
\put(80,35){\textcolor{black}{$\Gamma_{i+1}$}}
\put(145,35){\textcolor{black}{$\Gamma^p_{ i+1}$}}
\put(5,35){\textcolor{black}{$\Gamma^m_{ i+1}$}}
\end{picture}
 \caption{Illustration of the proxy surface \yzcmt{(dashed ellipse)} for compressing 
 $\mtx{A}_{i,i+1}$.   }
 \label{fig:offdiag_compression}
\end{figure}

While the actual matrix entries of \yzcmt{$\mtx{A}_{ij}$} are dependent on $\alpha$, the low rank factorization 
can be computed independent of $\alpha$ since {\color{black} the matrices need only be scaled by $\alpha$}.  As with 
the diagonal blocks, building the factorization of $\mtx{A}_{ij}$ directly is computationally prohibitive.
(The computational cost of the direct factorization is $O(N_iN_jk_{ij})$ where $k_{ij}$ is the numerical
rank of $\mtx{A}_{ij}$.)  Potential theory is again utilized to decrease the computational
cost.  
Consider an ellipse horizontally large enough to enclose $\Gamma_i$ 
and vertically shields $\Gamma_i$ from its top and bottom neighbor interface. A collection 
of $n_{\rm proxy}$ equivalent charges are evenly distributed on the ellipse in parameter space.  
Figure \ref{fig:offdiag_compression} illustrates
a proxy surface used for compressing $\mtx{A}_{i,i+1}$.  The interpolatory decomposition 
is applied to the matrix characterizing the interactions between the points on $\Gamma_i$ and
the proxy surface, $\mtx{A}^{\rm proxy}$.  The index vector $J_i$ and the $N_i\times k_{\rm proxy}$ interpolation
matrix \yzcmt{$\mtx{P}_{{\rm orig},ij}$} are returned.  Let \yzcmt{$J_{\rm orig} = J_i(1:k_{\rm proxy})$}. 

As with the diagonal block factorization, $k_{\rm proxy}$ is far from the optimal rank. To reduce the 
rank, we apply the interpolatory decomposition to \yzcmt{$[\mtx{A}_{0,ij} | \mtx{A}^p_{ij} | \mtx{A}^m_{ij}](J_{\rm orig},:)$}.
An $k_{\rm proxy}\times k_{\rm new}$ interpolation matrix $\mtx{P}_{{\rm new},ij}$ and index vector 
$J_{\rm new}$ are returned.  Then low rank factorization is complete. One factor can be used 
for all Bloch phases; \yzcmt{$\mtx{L}_{ij} = \mtx{P}_{{\rm orig},ij}\mtx{P}_{{\rm new},ij}$}.
The other factor
is simply a matrix evaluation; $\mtx{R}_{ij} = \mtx{A}_{ij}(J_{ij},:)$ where \yzcmt{$J_{ij}=J_{\rm orig}(J_{new})$}.  
It is important to note that the matrices $\mtx{A}_{0,ij}(J_{ij},:)$, {\color{black}$\mtx{A}^p_{ij}(J_{ij},:)$} 
and {\color{black}$\mtx{A}^m_{ij}(J_{ij},:)$} are computed once as they are independent of Bloch phase. Thus
constructing $\mtx{R}_{ij}$ \yzcmt{is formed simply by} matrix addition for each new Bloch phase.

{\color{black}\begin{remark}
 The rank of the factorizations of $\mtx{A}_{ij}$ will depend on the distance between
 the interfaces.   If the interfaces are space filling, the interaction between 
 the interfaces is not low rank.
\end{remark}
}

\subsection{The Bloch phase and incident angle dependence}
\label{sec:reuse}
Beyond the matrix $\mtx{A}$ and exploiting the sparsity of the 
other matrices in (\ref{eq:bigblock}), additional acceleration
can be gained for problems where the solution is desired for
multiple incident angles.

The matrix $\mtx{B}$ has {\color{black}block} entries (\ref{eq:Bmat}) that are independent of Bloch phase and thus need only be computed 
once.  This is also the case for $\mtx{V}$.  
The non-zero block matrices in $\mtx{C}$ \yzcmt{(\ref{eq:Cmat1}) and (\ref{eq:Cmat2}) are} dependent on
$\alpha$ but only as a constant multiple.  Thus the submatrices of $\mtx{C}$ can be precomputed and used
for all incident angles.  The same statement is true for \yzcmt{$\mtx{Q}$ and $\mtx{Z}$}. 

The only matrix that has entries that are dependent on incident angle is $\mtx{W}$.  In fact,
the matrix $\mtx{W}$ is only dependent on the Bloch phase $\alpha$.  Recall that 
the Bloch phase is defined as $\alpha = e^{id\omega_1 \cos \theta^{\rm inc}}$.  This means that 
for all incident angles that share a Bloch phase, there exist a representative 
angle $\hat{\theta}$ such that 
$\omega_1 \cos \theta^{\rm inc} = \hat{\theta} + \frac{2\pi m}{d}$ for \yzcmt{some} $m\in\mathbb{Z}$.  
Since the entries of $\mtx{W}$ involve 
$$e^{i\kappa_j x} = e^{i(\omega_1\cos \theta^{\rm inc}+\frac{2\pi j}{d})x}$$
for $j = -K,\ldots,K$.  Thus for angles with a shared
Bloch phase, the entries of $\mtx{W}$ are the same up to a shift in 
the index.  For example, suppose
that we know that 12 incident angles $\{\theta^{\rm inc}_1,\ldots,\theta^{\rm inc}_{12}\}$
share a Bloch phase and $\omega_1 \cos \theta^{\rm inc}_j = \hat{\theta} + \frac{2\pi (j-1)}{d}$
for $j = 1,\ldots, 12$.   The matrix $\mtx{W}$ is constructed so that it has entries with 
$\kappa_j$ indexed from $-K$ to $K+12$.  
This allows the singular value decomposition of $\mtx{S}$ to be 
used for all angles that share a Bloch phasewhich
results in substantial savings.  
To evaluate the solution using the resulting coefficients for the Rayleigh-Bloch expansion
above or below the unit cell, it is only necessary to use the terms that 
correspond to $-K,\ldots,K$ for that incident angle.

\subsection{Extensions}
\label{sec:extensions}
Many applications consider boundary value problem (\ref{eq:basic}) for
a collection of geometries where the variation is in a single interface 
or \high{wave number} in a layer.  The proposed direct solver can efficiently 
update an existing fast direct solver for these localized changes in the
geometry.

For example, if a user wants to replace $\Gamma_i$, {\color{black} only
the matrices corresponding to that interface need to be recomputed.  
This includes: the parts of the fast direct solver for $\mtx{A}$ 
corresponding to that block row and column, the corresponding block columns
of $\mtx{C}$ and the corresponding block rows of $\mtx{B}$.}  If the replaced 
layer is either the top or the bottom,
then sub-blocks in $\mtx{V}$ and $\mtx{Z}$ need to be updated as well.  
Independent of the interface
changed, the cost of creating a new direct 
solver is linear with respect to the number of discretization points on the new interface.  

If a user wants to change the wave number $\omega_i$ in $\Omega_i$ where $1<i<I$, there are two interfaces
affected $\Gamma_i$ and $\Gamma_{i+1}$.  The corresponding blocks rows and columns of the fast
direct solver of $\mtx{A}$ need to {\color{black} be} recomputed.  
In addition to updating those matrices, the corresponding blocks in $\mtx{B}$, $\mtx{C}$ and $\mtx{Q}$
need to {\color{black}be} updated.  If the wave number is changed in the top or bottom layer, then the 
corresponding blocks in $\mtx{Z}$ and $\mtx{V}$ need to be updated as well.  Again the cost of updating
the direct solver scales linearly with number of discretization points on the interfaces affected,
\yzcmt{i.e.}
the cost of the update is $O(N_i+N_{i+1})$.

\section{Numerical examples}
\label{sec:numerics}

This section illustrates the performance of the fast direct solver for 
 several \textcolor{black}{geometries with up to 11 layers, though the}
 solution technique can be applied to geometries of arbitrary number of layers.
Section \ref{sec:scaling} demonstrates the scaling of the fast direct solver for 
\textcolor{black}{ 3-layer and 9-layer geometries}
 where the wave number 
\textcolor{black}{alternates between $10$ and $10\sqrt{2}$ in the layers}.
  The ability for the solution technique to {\color{black}efficiently} solve (\ref{eq:basic}) for 
  hundreds of incident waves is {\color{black}demonstrated in section \ref{sec:angle_sweep}}.  
Section \ref{sec:perturbation} illustrates the performance of the solution technique {\color{black}
when the problem has: an interface geometry that is changed 
\yzcmt{or} a wave number that is changed in one of the layers.}


All {\color{black} the} geometries {\color{black} considered in this section} have period fixed at $d=1$.  
The vertical separation between the neighbor interfaces is roughly 1
\textcolor{black}{for the geometries considered in section \ref{sec:scaling} and \yzcmt{roughly} 1.5 for the rest of the experiments}. 
The interfaces are discretized via the Nystr\"om method with a 16-point Gaussian \yzcmt{composite} quadrature.  
The diagonal blocks
require specialized quadrature to handle the weakly singular kernels. {\color{black} For the experiments
presented in this section \yzcmt{generalized Gaussian quadrature \cite{gen_quad} is} utilized, but the fast direct solver is compatible with other specialized quadrature 
including Alpert \cite{alpert_quad}, Helsing \cite{helsing},
Kapur--Rokhlin  \cite{Kapur}, and QBX \cite{QBX}.  The geometries under consideration
involve both smooth interfaces and interfaces that have corners.}  
In order to achieve high accuracy without over discretizing, each 
corner is discretized with five levels of dyadic refinement.  Additionally, the 
integral operators are discretized in $\mathcal{L}^2$ \cite{Bremer1}. 
The artificial separation walls and proxy circles are discretized with 
parameter choices similar to those in \cite{Cho:15}.
\yzcmt{Specifically,}  
the left and right (vertical) artificial walls for each layer are discretized by $M_w=120$ \high{composite Guass Gauss-Legendre quadrature nodes}, 
the (horizontal) upper and lower walls of the unit cell are sampled at $M=60$ equispaced nodes 
and $P=160$ equispaced nodes are chosen on the proxy circle for each layer.  
For the wave numbers under consideration in these experiments, it is sufficient to truncate the 
Rayleigh-Bloch expansions at $K =20$.  

For all experiments, a HBS fast direct solver with tolerance $\epsilon = 10^{-12}$ was used to construct the approximation of
$\mtx{A}_0^{-1}$ in (\ref{eq:inv}).  The tolerance for all of the low rank factorizations was set to $10^{-12}$.
The singular value decomposition of $\mtx{S}$ was truncated at $\epsilon_{\rm Schur} = 10^{-13}$. 

All experiments were run on a dual 2.3 GHz Intel Xeon Processor E5-2695 v3 desktop workstation 
with 256 GB of RAM.  The code 
\yzcmt{is} implemented in MATLAB, apart from 
\yzcmt{the interpolatory decomposition, which uses Fortran}.  

The computational cost of the direct solution technique is broken into four parts:
\begin{itemize}
\item Precomputation I:  This consists of all computations for the \textit{fast linear algebra} that are independent of Bloch phase. This includes the fast application of $\mtx{A}_0^{-1}$, and the 
low rank factors $\mtx{L}_{ij}$ and $\mtx{R}_{ij}$ \yzcmt{needed} to make $\mtx{L}$ and $\mtx{R}$ as
presented in section \yzcmt{\ref{sec:fast_Ainv}}. The computational cost of this step is $O(N)$ where $N = \sum_{l=1}^I N_l$, \yzcmt{and} $N_l$ denotes
the number of discertization points on interface $l$.
 \item Precomputation II: This consists of the remainder of the precomputation that is independent of Bloch phase
 as presented in section \ref{sec:reuse}. The computational cost of this step is $O(N)$.
\item Precomputation III:  This consists of all the precomputation that can be used for incident angles that share a Bloch phase $\alpha$, including scaling matrices by $\alpha$, construction 
of the matrix $\mtx{W}$ as explained in section \ref{sec:reuse}, constructing the fast apply
of $\mtx{A}^{-1}$, evaluating the Schur complement matrix $\mtx{S}$ (\ref{eq:Schur}), and computing
the $\epsilon_{\rm Schur}$ SVD of $\mtx{S}$. The computational cost of this step is $O(N)$.
For a fixed number of discretization points on an interface, the computational cost is $O(I^3)$.
\item Solve: This consists of the application of the precomputed solver to the right hand side of (\ref{eq:bigblock}) via (\ref{eq:sigma}) and (\ref{eq:other}). 
The computational cost of the solve is $O(N)$.  For a fixed number of discretization points on an interface, the computational cost is $O(I^3)$.
\end{itemize}

The error is approximated via \yzcmt{a} flux error estimate as in \cite{Cho:15} which measures 
conservation of energy. This has been demonstrated to
agree with the relative error at any point in the domain.

\subsection{Scaling experiment}
\label{sec:scaling}
\textcolor{black}{This section illustrates the scaling of the fast direct solver for 
problems with \yzcmt{3-layers} and 9-layers corresponding to two and 
eight interface geometries. }

\begin{figure}[h]
\begin{picture}(200,200)(100,0)
\put(50,0){\includegraphics[trim={1 2cm 1 2cm},clip,scale=0.8]{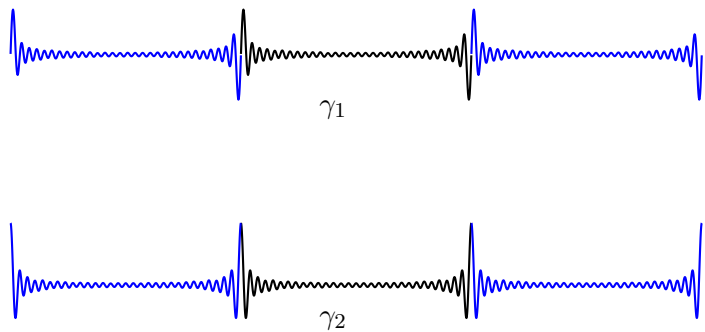}}
\put(210,105){$\gamma_1$ }
\put(210,25){$\gamma_2$ }
\end{picture}
\caption{
\textcolor{black}{Three periods of the interface geometries $\gamma_1$ and $\gamma_2$ as defined in equation (\ref{equ:fouriermap}). }
}
\label{fig:fouriermap_geometry}
\end{figure}

\yzcmt{For} the experiments \yzcmt{in} this section, the wave number in 
the layers remains fixed 
\textcolor{black}{(alternating between $10$ and $10\sqrt{2}$)} 
while the number of discretization points per layer increases.
 \textcolor{black}{The geometry consists of alternating the 
 following two interface geometries: $\gamma_1=(x_1(t),y_1(t))$ and 
 $\gamma_2=(x_2(t),y_2(t))$ \yzcmt{defined as}
\begin{equation}\label{equ:fouriermap}
\gamma_1:\,
 \begin{cases}
 x_1(t)=t-0.5\\
 y_1(t) = \frac{1}{60}\sum_{j=1}^{30}  a_j \sin(2\pi j t)\\
 \end{cases}
 \; \mbox{ and } \;
 \gamma_2:\,
 \begin{cases}
 x_2(t)=t-0.5\\
 y_2(t) = \frac{1}{60}\sum_{j=1}^{30}  b_j \cos(2\pi j t)\\
 \end{cases}
\end{equation}
for $t\in[0,1]$,
 where $\{a_j\}_{j=1}^{30}$ and $\{b_j\}_{j=1}^{30}$ are 
 random numbers in $[0,1)$ sorted in descending order.
Figure \ref{fig:fouriermap_geometry} illustrates the two 
interface geometries.
In each experiment, $\gamma_1$ and $\gamma_2$ are discretized with 
the same number of points $N_i$.  
The run time in seconds and flux error estimates are reported in Table \ref{tab:timing11}.
}

Each part of the 
solution technique scales linearly with respect to \textcolor{black}{$N_i$}.  
{\color{black}The factor of four increase in time for Precomputation I is 
expected since the cost scales linearly with the number of interfaces. 
Precomputation II should scale linearly with the number of layers and this
is observed with a factor three increase in the timings for this portion 
\yzcmt{of} the solver.  Precomputation III is expected to observe a factor nine increase
in the computational cost as this step scales cubically with the number of layers.
A factor of six is observed.  This is likely because the \yzcmt{problems under consideration are} sub-asymptotic in 
the number of layers. The same statement is true for the solve step.  
As 
expected the precomputation dominates the cost of the solver for both experiments.
 Precomputation parts I, II and III
account for approximately 
$90\%$, $3\%$, and $7\%$ of the precomputation time, respectively.
Thus the Bloch phase independent parts of the direct solver dominate
the computational cost.  }

\begin{table}[]
 \centering
 \begin{tabular}{|c|c|c|c|c|c|c|}\hline
  &$N_i$ & 1280 & 2560 & 5120 & 10240 & 20480 \\ \hline \hline
  \multirow{2}{*}{Precomp I} & 3-layer & 50 & 100 & 185 & 337 & 569 \\ \cline{2-7}
  & 9-layer &201 & 407& 768&1390& 2360 \\ \hline
  \multirow{2}{*}{Precomp II} & 3-layer &1.3 &2.4 & 4.5&8.2 &15.8\\ \cline{2-7}
  &9-layer& 3.8 & 7.6&13.4 & 26.1 & 51.4 \\ \hline
  \multirow{2}{*}{Precomp III}&3-layer & 2.5&5.9 & 10.6&21.4 & 41.9 \\ \cline{2-7}
&9-layer & 14.9 & 31.5 & 61.5 & 117.1 & 231.2 \\ \hline  
  \multirow{2}{*}{Solve} & 3-layer &0.1 &0.2 &0.8&1.6&3.4\\ \cline{2-7}
  &9-layer & 0.6& 2.5 & 4.9&13.5& 27.9\\ \hline
  \multirow{2}{*}{Flux error} & 3-layer &4.2e-5 & 6.9e-6 & 2.3e-8 &3.8e-10&4.5e-10\\ \cline{2-7}
  &9-layer & 2.1e-4 & 1.2e-5 & 1.5e-7 &4.6e-11 & 4.6e-10\\ \hline
  \end{tabular}
    \caption{\color{black}Time in seconds and flux error estimates for applying the direct solver to a \yzcmt{3-layer} and 
 9-layer geometry where the interfaces alternate between $\gamma_1$ and $\gamma_2$ defined in (\ref{equ:fouriermap}).
 $N_i$ denotes the number of discretization points for each boundary charge density on the interface.
    The wave number alternates between $10$ and $10\sqrt{2}$.
    }
    \label{tab:timing11}
  \end{table}

\subsection{Sweep over multiple incident angles}
\label{sec:angle_sweep}
\begin{figure}[h]
\begin{center}
\includegraphics[trim={1 0cm 0 0cm},clip,scale=0.72]{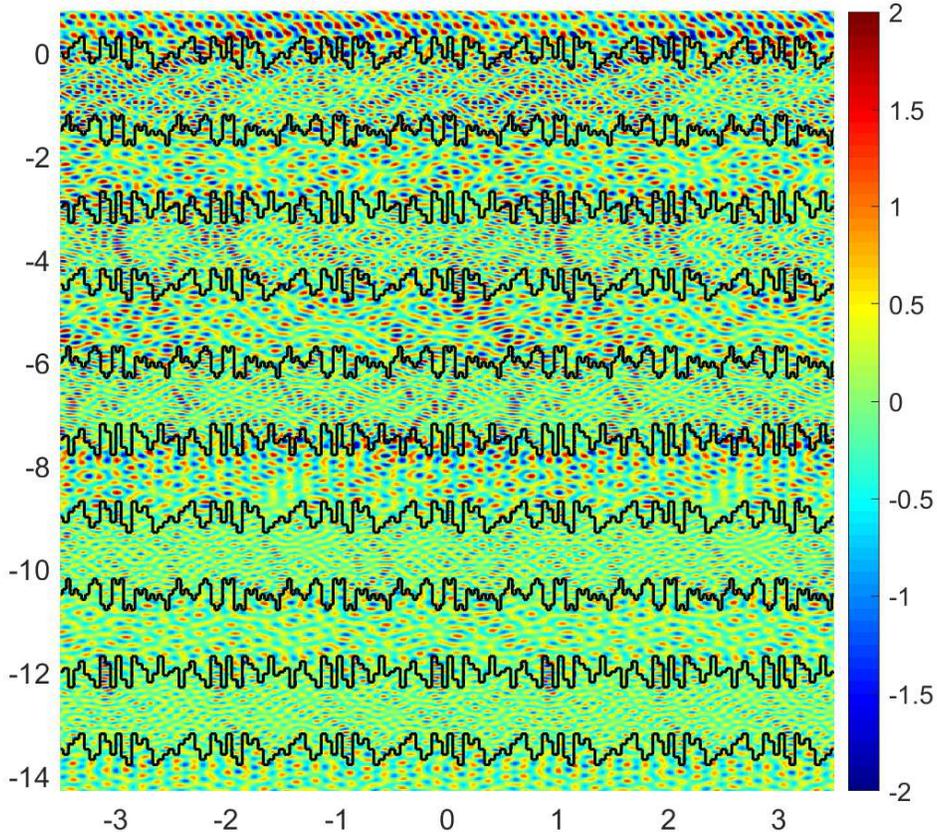}
\caption{Illustration of the real part of the total field of the solution to (\ref{eq:basic}) 
for a geometry with $10$ interfaces 
where the wave number alternates between $40$ and $40\sqrt{2}$.
The shown solution is for 
$\theta^{inc}=-0.845\pi$. 
 The total number of discretization 
points was set to $N =121,136$, resulting in a flux error estimate of $ 2.3e-8$. 
Seven periods in the geometry are shown. }
\label{fig:field_corner}
\end{center}
\end{figure}

\begin{figure}[h]
\begin{picture}(200,250)(100,0)
\put(70,0){\includegraphics[trim={0 0cm 0 0cm},clip,scale=0.5]{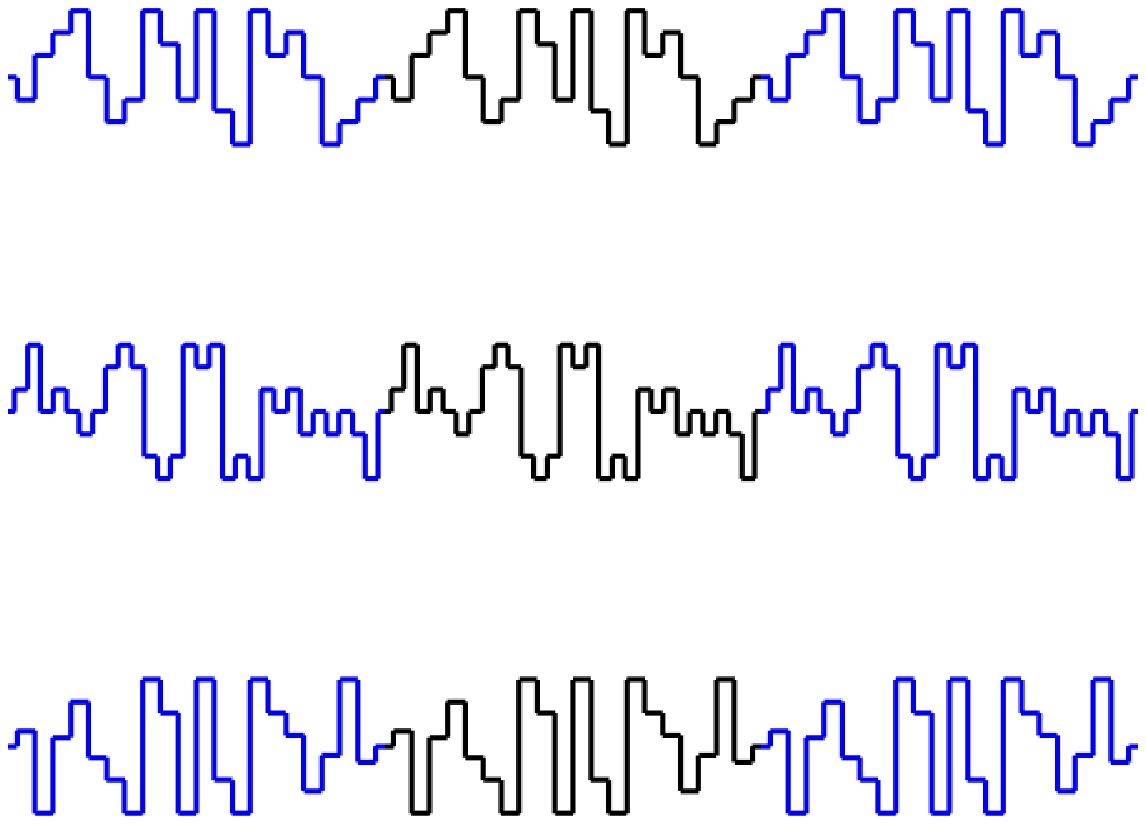}}
\put(70,105){`corner-1' geometry, 42 corners, 10,448 discretization points }
\put(70,55){`corner-2' geometry, 58 corners, 15,104 discretization points }
\put(70,5){`corner-3' geometry, 42 corners, 11,344 discretization points }
\end{picture}
\caption{The three different ``corner'' geometries in the 11-layer structure. 
Three periods are shown. 
See Figure \ref{fig:field_corner} for the full structure. }
\label{fig:corner_geometry}
\end{figure}

Many applications require solving (\ref{eq:basic}) for many incident angles
(as discussed in section \ref{sec:intro}).  In this setting, {\color{black}Precomputation} I and II \high{only needs to} 
be done once.  Precomputation III can be utilized for all incident angles 
that share a Bloch phase $\alpha$.  
\textcolor{black}{This section demonstrates the efficiency of the 
fast direct solver for handling scattering problems involving 
many incident angles. 
Specifically, we consider the geometry in Figure \ref{fig:field_corner}
which has eleven layers.  The interfaces consist of three different
corner geometries repeated in order.  
Each of the interfaces contains 40 to 50
right-angle corners. With the five levels of dyadic refinement into
each corner there are 10,000 to 15,000 discretization points per interface.
Figure \ref{fig:corner_geometry} provides \yzcmt{more details} about the
corner geometries including how many discretization 
points were used on each geometry.
The wave number in the layers alternates between $40$ and $40\sqrt{2}$. 
Equation (\ref{eq:basic}) was solved for $287$ incident angles between $[-0.89\pi, -0.11\pi]$.
This corresponds to 24 different 
Bloch phases.
 }
\textcolor{black}{
The average flux error estimate over the 287 incident angles 
is $2.4e-8$. 
Figure \ref{fig:field_corner} illustrates the real part of the total field for the 
incident angle $\theta^{inc}=-0.845\pi$.
The time for constructing and applying the fast direct solver for one incident angle 
is reported in the column labeled \textit{Original Problem} in Table \ref{tab:new}.  
Table \ref{tab:time_angle} reports the time for applying the proposed solution 
technique to the 287 boundary value problems using the same Precomputation I and 
II for all solves and exploiting the shared Bloch phase accelerations.  
}
Since the first two parts of the precomputation dominate the computational cost
for this geometry,
significant speed up over building the direct solver for each angle independently 
is observed.  For this problem, {\color{black}the proposed solution technique
is $100$ times }faster than building a fast direct solver
from scratch for each incident angle.  {\color{black} There is a $10$ 
times speed up
in the time for applying the solver for the multiple incident angles.  This 
results from the fact that angles that share a Bloch phase are processed together.}

  \begin{table}[]
    \centering
   \begin{tabular}{|c|l|l|l|l|l|}
    \hline
 $N_{total}$& Precomp I & Precomp II & Precomp III & Solve \\
     \hline  
121136 & 2369.3 & 32.3& 4517.4 & 482.2\\
           &            &        &(188.2 per Bloch phase) &( 1.7 per incident angle)\\
\hline
    \end{tabular}
    \caption{Time in seconds for solving 287 incident angles and 24 distinct Bloch phases
    on an 11-layer geometry shown in Figure \ref{fig:field_corner}. The incident angles are sampled from $[-0.89\pi, -0.11\pi]$. 
}
    \label{tab:time_angle}
    \end{table}

\subsection{Local change to the geometry}
\label{sec:perturbation}
{\color{black}This section illustrates the performance of the direct solver when there 
is a change in one layer of the geometry for the boundary value problem.  Either 
there is a change in an interface \yzcmt{geometry} or the wave number has been changed in a layer.}

\begin{figure}[h]
\centering
\begin{tabular}{cc}
 \includegraphics[trim={5cm .5cm 5cm .5cm},clip,width =4cm]{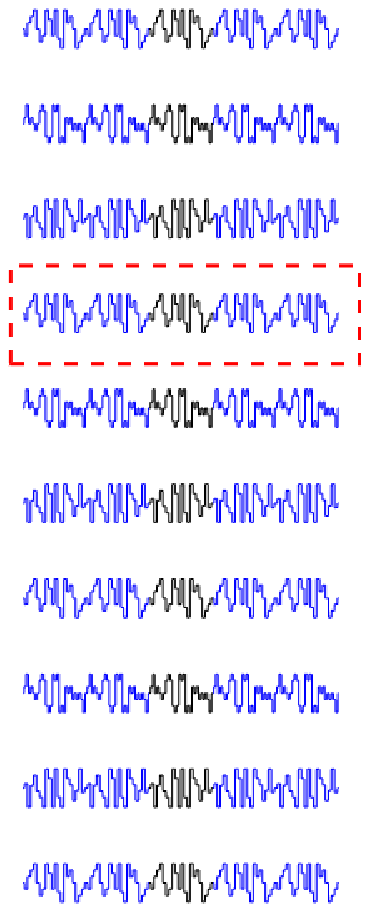}
 &  \includegraphics[trim={5cm .5cm 5cm .5cm},clip,width= 4cm]{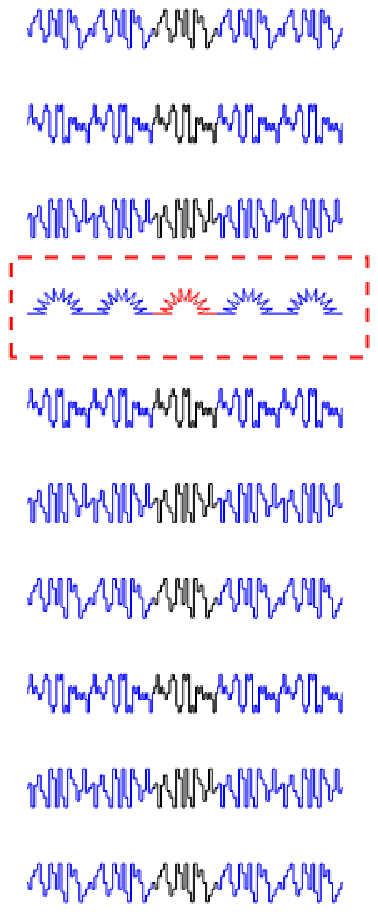}\\
 \yzcmt{(a)} & \yzcmt{(b)}
\end{tabular}
\caption{Illustration of 5 periods of (a) the original 11-layer structure and 
(b) the new structure obtained from replacing the fourth interface with a different geometry. 
The modified interface is in red box. 
}
\label{fig:newinnterface}
\end{figure}

We consider \textcolor{black}{an 11-layer geometry} 
\textcolor{black}{where the original set of interfaces are as illustrated in Figure \ref{fig:field_corner}}. 
{\color{black}As in the last section, the wave number alternates between $40$ and $40\sqrt{2}$ 
in the layers of the original geometry.}
\textcolor{black}{The incident angle for these experiments is fixed at $\theta^{inc}=-\frac{\pi}{5}$.}
In the first experiment, the fourth interface in \textcolor{black}{the original} geometry is changed
to the ``hedgehog'' interface.  Figure \ref{fig:newinnterface} illustrates the original and new geometries.  
 The hedgehog geometry \yzcmt{consists} of 17 sharp corners and cannot
be \yzcmt{written} as the graph of a function defined on the $x$-axis.  The number of discretization
points \yzcmt{on the new interface} needed to maintain the same accuracy as the original problem is
 $N_4=14,496$.  
In the second experiment, the wave number for the second layer is changed from 
$40\sqrt{2}$ to $30$ but the interfaces are kept fixed with the geometry illustrated in 
Figure \ref{fig:newinnterface}(a).
As presented in section \ref{sec:reuse}, only 
a small number of the matrices in each step of the precomputation need to be 
recomputed.  

{\color{black}Table \ref{tab:new} reports time in seconds for building a direct solver
from scratch for the original problem, updating the solver when the 
fourth interface is replaced and updating the direct solver when there
is a change in wave number in the second layer.  
The parts of Precomputation I and II that are needed in the updating scheme
are smaller than building them from scratch.  Precomputation I is approximately
twice as expensive for the problem with the changed wave number because it requires
updating matrices for two interfaces while the changed interface problem only 
involves one interface.
  Precomputation III needs to be redone for the new wave number problem.  
  This is why it is nearly as expensive as
for the original problem.  There is slight savings because the ranks related
to \yzcmt{the} replaced wave number \yzcmt{are} lower since the new wave number is smaller than 
the original. 
\yzcmt{The cost savings for updating the solver} is greater for the changed
interface problem since more of the solver from the original problem 
can be re-used.}

 When an application 
requires many solves per new geometry and many new geometries need to be 
considered, the speed gains over building a solver from scratch for each 
new geometry will be significant.  

 \begin{table}[h]
\centering
\begin{tabular}{|c|c|c|c|
}
\hline
&  Original problem &Replace interface $\Gamma_4$ & Change wave number $\omega_2= 30$ \\ \hline
$N_{total}$&121,136 &   125,184  &121,136 \\ \hline
Precomp I &2369.3& 237 &433  \\ \hline
Precomp II &32.3&11.7  & 3.5 \\ \hline
Precomp III  &174.1& 41.7 & 109.2  \\ \hline
Solve &18.8&15.7& 13.6 \\ 
\hline 
\end{tabular}
\caption{Time in seconds for constructing and applying the fast direct solver 
\yzcmt{to an} 11-layer geometry (first column), a geometry that 
has the fourth interface \yzcmt{changed} (second column) and the wave number for the second layer 
changed from $40\sqrt{2}$ to $30$ (third column). $N_{total}$ is the number of discretization points on the 
interfaces in the unit cell.  }
\label{tab:new}
\end{table}%

\section{Conclusion}
\label{sec:summary}
This paper presents a fast direct solution technique for multilayered media
quasi-periodic scattering problem.  For low frequency problems, the computational
cost of the direct solver scales linearly with the number of discretization points.
The bulk of the precomputation can be used for all solves independent of incident 
angle and Bloch phase $\alpha$.  For a problem where over 200 hundred
incident angles are considered, the proposed fast direct solver is 100 times faster
than building a direct solver from scratch.

An additional benefit of this solution technique is that modifications
in \yzcmt{the wave number of a layer} or an interface \yzcmt{geometry} result in only having to update 
the matrices corresponding to that layer or interface.  The cost of updating
the precomputation parts scales linearly with the number of points on that
interface.  For a problem with a changed interface, the constant associated
with the linear scaling is very small for the precomputation (relative to 
building a new direct solver from scratch).  In optimal 
design and inverse scattering applications where 
the geometry will be changed many times and for each geometry many solves
are required, the fast direct solver will have significant savings.  

Two dimensional geometries have to be  
complex in order to justify the need for the fast direct solver.  
For three dimensional problems, a fast direct solver will be necessary for
\high{most geometries of interest in applications.}
The extension to three dimensional problems is not trivial but the work
presented in this paper provides the foundations for that work.

\section*{Acknowledgements}
The authors thank Alex Barnett for the fruitful conversations.  
The work of A. Gillman is supported by the National Science Foundation (DMS-1522631).  

\appendix
\section{The construction of $\mtx{S}_2$}
\label{app:1}

This section presents an efficient technique for constructing the tridiagonal matrix 
$\mtx{S}_2 = \mtx{I} +\mtx{R}\mtx{A}_0^{-1}\mtx{L}$.  

For simplicity of presentation, \yzcmt{let the blocks of $\mtx{S}_2$ be denoted as follows}
\begin{itemize}
 \item[] $\mtx{X}_i$ for $1\leq i\leq I$ denotes the diagonal blocks,
 \item[] $\mtx{Y}_i$ for $2\leq i \leq I$ denotes the lower diagonal blocks, and
 \item[] $\mtx{Z}_i$ for $1\leq i\leq I-1$ denotes the upper diagonal blocks.
\end{itemize}

The diagonal blocks are given by 
$$
\mtx{X}_1=
\begin{bmatrix}
\mtx{I}+{\color{black}\mtx{R}^{pm}_{11}} \mtx{A}_{0,11}^{-1} {\color{black}\mtx{L}^{pm}_{11}} & {\color{black}\mtx{R}^{pm}_{11}} \mtx{A}_{0,11}^{-1} \mtx{L}_{12}\\
\mtx{0} & \mtx{I}
\end{bmatrix}	,
$$
$$
\mtx{X}_I=
\begin{bmatrix}
\mtx{I} & \mtx{0}\\
{\color{black}\mtx{R}^{pm}_{II}} \mtx{A}_{0,II}^{-1} \mtx{L}_{I,I-1} & \mtx{I} +{\color{black}\mtx{R}^{pm}_{II}} \mtx{A}_{0,II}^{-1} {\color{black}\mtx{L}^{pm}_{II}}
\end{bmatrix}	,
$$
and, for $2\leq i\leq (I-1)$,  
$$
\mtx{X}_{i}=
\begin{bmatrix}
\mtx{I} & \mtx{0} & \mtx{0}\\
{\color{black}\mtx{R}^{pm}_{ii}} \mtx{A}_{0,ii}^{-1} \mtx{L}_{i,i-1} 
	& \mtx{I} + {\color{black}\mtx{R}^{pm}_{ii}} \mtx{A}_{0,ii}^{-1} {\color{black}\mtx{L}^{pm}_{ii}}
	&   {\color{black}\mtx{R}^{pm}_{ii}} \mtx{A}_{0,ii}^{-1} \mtx{L}_{i,i+1}\\
\mtx{0} & \mtx{0} &\mtx{I}
\end{bmatrix}	.
$$
The lower diagonal blocks are given by 
$$
\mtx{Y}_2=
\begin{bmatrix}
\mtx{R}_{21} \mtx{A}_{0,11}^{-1} {\color{black}\mtx{L}^{pm}_{11}} & 	\mtx{R}_{21} \mtx{A}_{0,11}^{-1} \mtx{L}_{12}\\
\mtx{0} & \mtx{0}
\end{bmatrix},
$$
and, for $3\leq i\leq I$,
$$
\mtx{Y}_i=
\begin{bmatrix}
\mtx{R}_{i,i-1} \mtx{A}_{0,(i-1)(i-1)}^{-1} \mtx{L}_{i-1,i-2} 
	& \mtx{R}_{i,i-1} \mtx{A}_{0,(i-1)(i-1)}^{-1} {\color{black}\mtx{L}^{pm}_{(i-1)(i-1)}}
	&  \mtx{R}_{i,i-1} \mtx{A}_{0,(i-1)(i-1)}^{-1} \mtx{L}_{i-1,i}\\
\mtx{0} & \mtx{0} &\mtx{0}
\end{bmatrix}	.
$$ 
Finally the upper diagonal blocks are defined by 

$$
\mtx{Z}_i=
\begin{bmatrix}
\mtx{0} & \mtx{0} & \mtx{0} \\
\mtx{R}_{i,i+1} \mtx{A}_{0,(i+1)(i+1)}^{-1} \mtx{L}_{i+1,i}
	&\mtx{R}_{i,i+1} \mtx{A}_{0,(i+1)(i+1)}^{-1} {\color{black}\mtx{L}^{pm}_{(i+1)(i+1)}}
	& \mtx{R}_{i,i+1} \mtx{A}_{0,(i+1)(i+1)}^{-1}  \mtx{L}_{i+1,i+2}
\end{bmatrix}	
$$
for $1\leq i\leq (I-2)$, 
and 
 $$
\mtx{Z}_{I-1}=	
\begin{bmatrix}
\mtx{0} & \mtx{0}\\
\mtx{R}_{I-1,I} \mtx{A}_{0,II}^{-1} \mtx{L}_{I,I-1} 
	& \mtx{R}_{I-1,I} \mtx{A}_{0,II}^{-1} {\color{black}\mtx{L}^{pm}_{ II}}
\end{bmatrix}. 
$$

The matrix $\mtx{S}_2$ can be inverted via a block variant of the Thomas algorithm.  
Let the sum of the ranks of the low-rank approximations be defined as $N_{i}^{block}=k_{pm,ii}+k_{i,i-1}+k_{i,i+1}$ for $2\leq i \leq I-1$, $N_1^{block}=k_{pm,11}+k_{1,2}$ and $N_I^{block}=k_{pm,II}+k_{I,I-1}$.
The diagonal block $\mtx{X}_i$ is of size $N_{i}^{block} \times N_{i}^{block}$.
The upper diagonal block $\mtx{Z}_i$ has size $N_{i}^{block} \times N_{i+1}^{block}$.
The lower diagonal block $\mtx{Y}_i$ has size $N_{i}^{block} \times N_{i-1}^{block}$.

\yzcmt{
For the tested geometries and wave numbers, $N_{i}^{block}$ is only several hundreds and 
the diagonal blocks can be 
inverted rapidly via dense linear algebra.}
If all of the blocks are of similar size $N_{i}^{block}\approx N^{block}$, 
then the cost \yzcmt{of} inverting $\mtx{S}_2$ via the block Thomas algorithm is 
\yzcmt{$\mathcal{O}([N^{block}]^3I)$,
which is linear with respect to the number of interfaces.}


\begin{thebibliography}{10}


\bibitem{Abramowitz-Stegun}
M.~Abramowitz and I.~Stegun, editors.
\newblock {\em Handbook of Mathematical Functions}.
\newblock Dover, New York, 1964.

\bibitem{2019_dielectric}
D. A. Alessi, H. T. Nguyen, J. A. Britten, P. A. Rosso, and C. Haefner.
\newblock Low-dispersion low-loss dielectric gratings for efficient ultrafast laser pulse compression at high average powers.
\newblock {\em Optics \& Laser Technology}, 117: 239--243, 2019.



\bibitem{alpert_quad}
B.~Alpert.
\newblock Hybrid gauss-trapezoidal quadrature rules.
\newblock {\em SIAM Journal on Scientific Computing}, 20(5):1551--1584, 1999.


\bibitem{Arens_thesis}
T. Arens. 
\newblock Scattering by biperiodic layered media: The integral equation ap-
proach. 
\newblock Habilitation thesis, Karlsruhe, 2010.


\bibitem{2006_Arens}
T. Arens, S. N. Chandler-Wilde, and J. A. DeSanto. 
\newblock On integral equation and
least squares methods for scattering by diffraction gratings. 
\newblock {\em Computer Physics Communications}, 1:1010--42, 2006.

\bibitem{2010_solar1}
H.~A. Atwater and A.~Polman./
\newblock Plasmonics for improved photovoltaic devices.
\newblock {\em Nature Materials}, 9:205-- 213, 2010.



\bibitem{pollution}
I.~M. Babuska and S.~A. Sauter.
\newblock Is the pollution effect of the {FEM} avoidable for the {H}elmholtz
  equation considering high wave numbers?
\newblock {\em SIAM Journal of Numerical Analysis}, 34(6):2392--2423, 1997.

\bibitem{timo}
A.~Barnett and T.~ Betcke.
\newblock Stability and convergence of the method of fundamental solutions for Helmholtz problems on analytic domains.
\newblock {\em Journal of Computational Physics}, 227(14): 7003--7026, 2008.

\bibitem{Barnett-Greengard}
A.~Barnett and L.~Greengard.
\newblock A new integral representation for quasi-periodic fields and its
  application to two-dimensional band structure calculations.
\newblock {\em Journal of Computational Physics}, 229:6898--6914, 2010.

\bibitem{2004_laser}
C.~Barty, M.~Key, J.~Britten, R.~Beach, G.~Beer, C.~Brown, S.~Bryan, J.~Caird,
  T.~Carlson, J.~Crane, J.~Dawson, A.~Erlandson, D.~Fittinghoff, M.~Hermann,
  C.~Hoaglan, A.~Iyer, L.~J. II, I.~Jovanovic, A.~Komashko, O.~Landen, Z.~Liao,
  W.~Molander, S.~Mitchell, E.~Moses, N.~Nielsen, H.-H. Nguyen, J.~Nissen,
  S.~Payne, D.~Pennington, L.~Risinger, M.~Rushford, K.~Skulina, M.~Spaeth,
  B.~Stuart, G.~Tietbohl, and B.~Wattellier.
\newblock An overview of llnl high-energy short-pulse technology for advanced
  radiography of laser fusion experiments.
\newblock {\em Nuclear Fusion}, 44(12):S266, 2004.

\bibitem{2010_borm_book}
S.~B{\"o}rm.
\newblock {\em Efficient numerical methods for non-local operators}, volume~14
  of {\em EMS Tracts in Mathematics}.
\newblock European Mathematical Society (EMS), Z\"urich, 2010.

\bibitem{2004_borm_hackbusch}
S.~B{\"o}rm and W.~Hackbusch.
\newblock Approximation of boundary element operators by adaptive
  {$\mathcal{H}^2$}-matrices.
\newblock In {\em Foundations of computational mathematics: {M}inneapolis,
  2002}, volume 312 of {\em London Math. Soc. Lecture Note Ser.}, pages 58--75.
  Cambridge Univ. Press, Cambridge, 2004.

\bibitem{Bremer1}
J.~Bremer.
\newblock On the {N}ystr\"om discretization of integral operators on planar
  domains with corners.
\newblock {\em Applied and Computational Harmonic Analysis}, 32:45--64, 2012.

\bibitem{2013_3DBIE}
J.~Bremer, A.~Gillman, and P.~Martinsson.
\newblock A high-order accurate accelerated direct solver for acoustic
  scattering from surfaces.
\newblock {\em BIT Numerical Mathematics}, 55:141--170, 2015.




\bibitem{2017_Bruno}
O.~P. Bruno and A.~G. Fernandez-Lado.
\newblock Rapidly convergent quasi-periodic green functions for scattering by
  arrays of cylinders{\textemdash}including wood anomalies.
\newblock {\em Proceedings of the Royal Society of London A: Mathematical,
  Physical and Engineering Sciences}, 473(2199), 2017.

\bibitem{2019_Bruno}
O. P. Bruno and M. C. Haslam. 
\newblock Efficient high-order evaluation of scattering by
periodic surfaces: deep gratings, high frequencies, and glancing incidences. 
\newblock {\em Journal of the Optical Society of America A}, 26(3):658--668, 2009.


  
\bibitem{lowrank}
H.~Cheng, Z.~Gimbutas, P.~Martinsson, and V.~Rokhlin.
\newblock On the compression of low rank matrices.
\newblock {\em SIAM Journal of Scientific Computing}, 26(4):1389--1404, 2005.

\bibitem{Cho:18}
M.~Cho.
\newblock Spectrally-accurate numerical method for acoustic scattering from
  doubly-periodic 3d multilayered media.
\newblock {\em Journal of Computational Physics}, 393: 46--58, 2019. 
  
\bibitem{Cho:15}
M.~Cho and A.~Barnett.
\newblock Robust fast direct integral equation solver for quasi-periodic
  scattering problems with a large number of layers.
\newblock {\em Optics Express}, 23(2):1775--1799, 2015.

\bibitem{RCWA4}
T.~H. Chou, K.~Y. Cheng, T.~L. Chang, C.~J. Ting, H.~C. Hsu, C.~J. Wu, J.~H. Tsai, and T.~Y. Huang.
\newblock Fabrication of antireflection structures on TCO film for reflective liquid crystal display.
\newblock {\em Microelectronic Engineering}, 86(4), 628--631, 2009.



\bibitem{coltonkress}
D.~Colton and R.~Kress.
\newblock {\em Inverse acoustic and electromagnetic scattering theory},
  volume~93 of {\em Applied Mathematical Sciences}.
\newblock Springer-Verlag, Berlin, second edition, 1998.

\bibitem{Darbe2019SimulationAP}
S.~Darbe, M.~D. Escarra, E.~C. Warmann and H.~A. Atwater.
\newblock Simulation and partial prototyping of an eight‐junction holographic spectrum‐splitting photovoltaic module.
\newblock {\em Energy Science \& Engineering}, 7(6), 2019. Web. doi:10.1002/ese3.445.

\bibitem{Miller}
V. Ganapati, O. D. Miller and E. Yablonovitch.
\newblock Light Trapping Textures Designed by Electromagnetic Optimization for Subwavelength Thick Solar Cells.
\newblock {\em IEEE Journal of Photovoltaics}, 4(1):175--182, 2014.


\bibitem{gmesh}
C.~Geuzaine and J.-F. Remacle.
\newblock Gmsh: a three-dimensional finite element mesh generator with built-in
  pre- and post-processing facilities.
\newblock {\em International Journal for Numerical Methods in Engineering},
  79(11):1309--1331, 2009.


  \bibitem{2013_Gillman}
  A. Gillman and A. Barnett.
  \newblock A fast direct solver for quasiperiodic scattering problems.
  \newblock {\em Journal of Computational Physics} 248:309--322, 2013.
  
  
  
\bibitem{2012_martinsson_FDS_survey}
A.~Gillman, P.~Young, and P.~Martinsson.
\newblock A direct solver {$O(N)$} complexity for integral equations on
  one-dimensional domains.
\newblock {\em Frontiers of Mathematics in China}, 7:217--247, 2012.

\bibitem{golub}
G.~H. Golub and C.~F. Van~Loan.
\newblock {\em Matrix computations}.
\newblock Johns Hopkins Studies in the Mathematical Sciences. Johns Hopkins
  University Press, Baltimore, MD, third edition, 1996.

\bibitem{2014_periodic_HO}
L.~Greengard, K.~Ho, and J.-Y. Lee.
\newblock A fast direct solver for scattering from periodic structures with
  multiple material interfaces in two dimensions.
\newblock {\em Journal of Computational Physics}, 258:738--751., 2014.

\bibitem{gu1996}
M.~Gu and S.~C. Eisenstat.
\newblock Efficient algorithms for computing a strong rank-revealing {QR}
  factorization.
\newblock {\em SIAM Journal on Scientific Computing}, 17(4):848--869, 1996.


\bibitem{RCWA1}
K.~Han, and C.~H. Chang.
\newblock Numerical Modeling of Sub-Wavelength Anti-Reflective Structures for Solar Module Applications.
\newblock {\em Nanomaterials}, 4(1): 87--128, 2014.

\bibitem{gen_quad}
S.~Hao, A.~H. Barnett, P.~G. Martinsson, and P.~Young.
\newblock High-order accurate nystrom discretization of integral equations with
  weakly singular kernels on smooth curves in the plane.
\newblock {\em Advances in Computational Mathematics}, 40:245--272, 2013.

\bibitem{helsing}
J.~Helsing and R.~Ojala.
\newblock Corner singularities for elliptic problems: integral equations,
  graded meshes, quadrature, and compressed inverse preconditioning.
\newblock {\em Journal of Computational Physics}, 227:8820--8840, 2008.

\bibitem{2012_ho_greengard_fastdirect}
K.~Ho and L.~Greengard.
\newblock A fast direct solver for structured linear systems by recursive
  skeletonization.
\newblock {\em SIAM Journal of Scientific Computing}, 34(5):2507--2532, 2012.

\bibitem{2014_HIF}
K.~Ho and L.~Ying.
\newblock Hierarchical interpolative factorization for elliptic operators:
  Integral equations.
\newblock {\em Communications on Pure and Applied Mathematics},
  69(7):1314--1353, 2015.
  
  
\bibitem{2002_Wilde}
K. V. Horoshenkov and S. N. Chandler-Wilde. 
\newblock Efficient calculation of two-
dimensional periodic and waveguide acoustic Green’s functions. 
\newblock {\em Journal of the Acoustical Society of America}, 111:1610--1622, 2002.  

\bibitem{Hughes}
T.~Hughes.
\newblock {\em The Finite Element Method: Linear Static and Dynamic Finite
  Element Analysis}.
\newblock Dover Publications, 2000.

\bibitem{2010_chirped}
G.~A. Kalinchenko and A.~M. Lerer.
\newblock Wideband all-dielectric diffraction grating on chirped mirror.
\newblock {\em Journal of Lightwave Technology}, 28:2743--2749, 2010.

\bibitem{Kapur}
S.~Kapur and V.~Rokhlin.
\newblock High-order corrected trapezoidal quadrature rules for singular
  functions.
\newblock {\em SIAM Journal of Numerical Analysis}, 34(4):1331--1356, 1997.

\bibitem{2010_solar2}
M.~D. Kelzenberg, S.~W. Boettcher, J.~A. Petykiewicz, D.~B. Turner-Evans, M.~C.
  Putnam, E.~L. Warren, J.~M. Spurgeon, R.~M. Briggs, N.~S. Lewis, and H.~A.
  Atwater.
\newblock Enhanced absorption and carrier collection in si wire arrays for
  photovoltaic applications.
\newblock {\em Nature Materials}, 9:239--244, 2010.

\bibitem{QBX}
A.~{Kl\"ockner}, A.~Barnett, L.~Greengard, and M.~O'Neil.
\newblock Quadrature by expansion: A new method for the evaluation of layer
  potentials.
\newblock {\em Journal of Computational Physics}, 252:332 -- 349, 2013.

\bibitem{2003_pml}
D.~Komatitsch and J.~Tromp.
\newblock A perfectly matched layer absorbing boundary condition for the
  second-order seismic wave equation.
\newblock {\em Geophysical Journal International}, 154(1):146--153, 2003.


\bibitem{1978_kress}
A.~Kress and G.~F. Roach.
\newblock Transmission problems for the Helmholtz equation.
\newblock {\em Journal of Mathematical Physics}, 19(6): 1433--1437, 1978.



\bibitem{1996_RCWA_li2}
L.~Li.
\newblock Formulation and comparison of two recursive matrix algorithms for
  modeling layered diffraction gratings.
\newblock {\em Journal of the Optical Society of America A}, 13:1024--1035,
  1996.

\bibitem{1996_RCWA_li1}
L.~Li.
\newblock Use of fourier series in the analysis of discontinuous periodic
  structures.
\newblock {\em Journal of the Optical Society of America A}, 13:1870--1876,
  1996.
  
  
\bibitem{2007_linton}  
C.~M. Linton and I.~Thompson.
\newblock Resonant effects in scattering by periodic arrays.
\newblock {\em Wave Motion}, 44:165--175, 2007.

\bibitem{2016_periodic_stokes}
G.~Marple, A.~Barnett, A.~Gillman, and S.~Veerapaneni.
\newblock A fast algorithm for simulating multiphase flows through periodic
  geometries of arbitrary shape.
\newblock {\em SIAM Journal of Scientific Computing}, 38(5):B740--B772, 2016.





\bibitem{1981_RCWA}
M.~G. Moharam and T.~G. Gaylord.
\newblock Rigorous coupled-wave analysis of planar-grating diffraction.
\newblock {\em Journal of the Optical Society of America}, 71:811--818, 1981.



\bibitem{2008_Nicholas}
M.~J. Nicholas.
\newblock A higher order numerical method for 3-D doubly periodic electromagnetic scattering problems.
\newblock Communications in Mathematical Sciences, 6: 669--694, 2008.

\bibitem{1995_diegrating}
M.~D. Perry, R.~D. Boyd, J.~A. Britten, D.~Decker, B.~W. Shore, C.~Shannon, and
  E.~Shults.
\newblock High-efficiency multilayer dielectric diffraction gratings.
\newblock {\em Optics Letters}, 20:940--942, 1995.

\bibitem{rokh83}
V.~Rokhlin.
\newblock Solution of acoustic scattering problems by means of second kind
  integral equations.
\newblock {\em Wave Motion}, 5:257--272, 1983.

\bibitem{Sergeant:10}
N.~Sergeant, M.~Agrawal, and P.~Peumans.
\newblock High performance solar-selective absorbers using coated
  sub-wavelength gratings.
\newblock {\em Optics Express}, 18(6):5525--5540, 2010.


\bibitem{2007_shiv_sheng}
Z.~Sheng, P.~Dewilde, and S.~Chandrasekaran.
\newblock Algorithms to solve hierarchically semi-separable systems.
\newblock In {\em System theory, the Schur algorithm and multidimensional
  analysis}, volume 176 of {\em Operator Theory: Advances and Applications},
  pages 255--294. Birkh\"auser, Basel, 2007.

\bibitem{RCWA3}
C.~J. Ting, C.~F. Chen, C.~P. Chou. 
\newblock Antireflection subwavelength structures analyzed by using the finite difference time domain method. 
\newblock {\em Optik}, 120:814--817, 2009.
  
  \bibitem{RCWA2}
H.~Y. Tsai. 
\newblock Finite difference time domain analysis of three-dimensional sub-wavelength structured arrays.
\newblock {\em Japanese Journal of Applied Physics}, 47:5007--5009, 2008. 
  
  
\bibitem{2010_xia}
J.~Xia, S.~Chandrasekaran, M.~Gu, and X.~Li.
\newblock Fast algorithms for hierarchically semiseparable matrices.
\newblock {\em Numerical Linear Algebra with Applications}, 17(6):953--976,
  2010.

\bibitem{2009_xia_superfast}
J.~Xia, S.~Chandrasekaran, M.~Gu, and X.~S. Li.
\newblock Superfast multifrontal method for large structured linear systems of
  equations.
\newblock {\em SIAM Journal on Matrix Analysis and Applications},
  31(3):1382--1411, 2009.
  
  \bibitem{BBStarling1994}
  A-S.~Bonnet-Bendhia and F.~Starling.
  \newblock Guided waves by electromagnetic gratings and non-uniqueness examples for the diffraction problem.
  \newblock{\em Mathematical Methods in the Applied Sciences}, 
  17(5):305-338, 1994.

  \bibitem{Yip_woodburystability}
E.L.~Yip.
\newblock A note on the stability of solving a rank-p modification of a linear system by the
Sherman-Morrison-Woodbury formula.
\newblock{\em SIAM Journal on Scientific and Statistical Computing},
7(2):507-513, 1986.
\end{thebibliography}
\end{document}